\documentclass[]{report}
\usepackage{amssymb}
\usepackage{titling}
\usepackage{amsmath}
\usepackage{geometry}
\usepackage{graphicx}
\usepackage{subfig}
\usepackage{fancyhdr}
\usepackage{lastpage}
\usepackage{lipsum}
\usepackage{fancyhdr}
\usepackage{lastpage}
\usepackage{lipsum}
\usepackage{ragged2e}
\usepackage{amsmath}
\usepackage[ruled,vlined]{algorithm2e}
\usepackage{algpseudocode}
\usepackage{comment}

\newtheorem{theorem}{Theorem}
\usepackage{hyperref}
\hypersetup{colorlinks,
	citecolor=blue
}
\usepackage{listings}
\usepackage{color}

\definecolor{dkgreen}{rgb}{0,0.6,0}
\definecolor{gray}{rgb}{0.5,0.5,0.5}
\definecolor{mauve}{rgb}{0.58,0,0.82}

\title{Numerical Analysis for Real-time Nonlinear
	Model Predictive Control (NMPC) of Ethanol
	Steam Reformers}
\author{Robert Joseph\\ Department of Mathematical and Statistical Sciences}

\begin{document}


\begin{titlepage}
    \begin{center}
        \vspace*{1cm}
            
        \huge
        \textbf{Numerical Analysis for Real-time Nonlinear
	Model Predictive Control of Ethanol
	Steam Reformers}
            
        \vspace{0.5cm}
        \LARGE
        An undergraduate thesis presented by            
        \vspace{0.5cm}  
        
        {Robert Joseph George}
        \vspace{0.5cm}
        
        Advised by: Xinwei Yu\\
        Supervised by: Paul Buckingham
        \vspace{1cm}

        \normalsize{\hspace{12pt} \textbf{Abstract}}
        
        \justify 
        {The utilization of renewable energy technologies, particularly hydrogen, has seen a boom in interest and has spread throughout the world. Ethanol steam reformation is one of the primary methods capable of producing hydrogen efficiently and reliably. This paper provides an in-depth study of the reformulated system, both theoretically and numerically, as well as a plan to explore the possibility of converting the system into its conservation form. Lastly, we offer an overview of several numerical approaches for solving the general first-order quasi-linear hyperbolic equation to the particular model for ethanol steam reforming (ESR). We conclude by presenting some results that would enable the usage of these ODE/PDE solvers to be used in non-linear model predictive control (NMPC) algorithms and discuss the limitations of our approach and directions for future work.}

        \vspace{1cm}
        \centering
        \normalsize{Submitted in partial fulfillment of the Honors requirements for the degree \\ of Bachelor of Honors in Applied Mathematics and Computer Science.}
            
        \vspace{0.8cm}
	\includegraphics[width=0.7\linewidth]{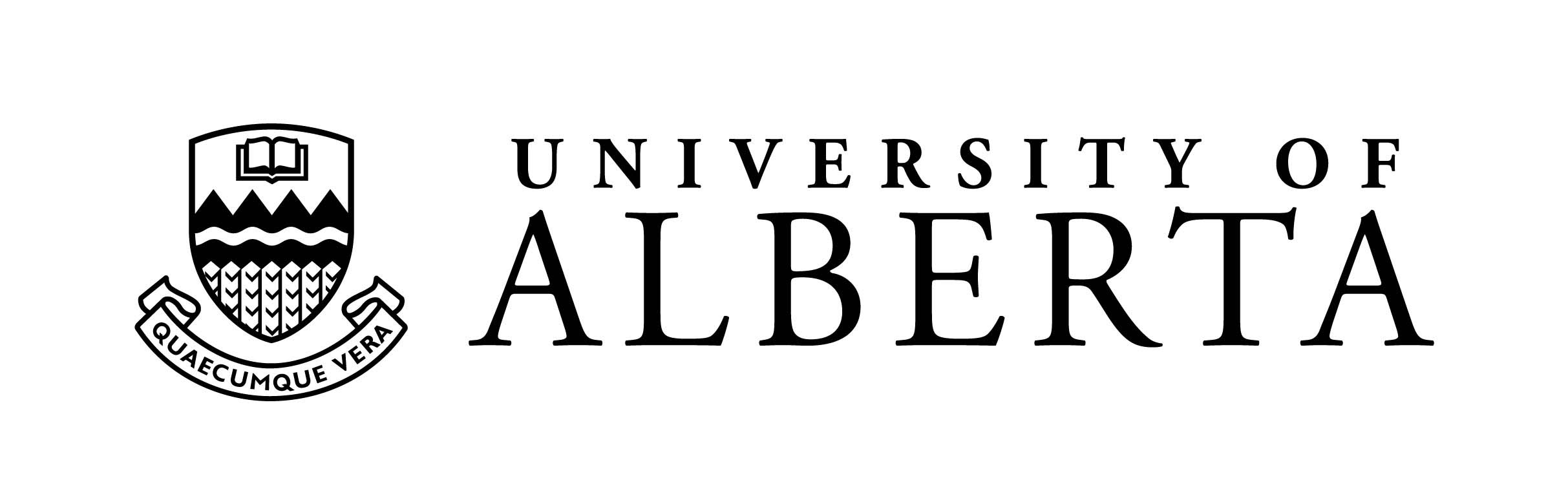}

        \Large
        Department of Mathematics and Statistics \\
        \hspace{12pt} Edmonton, Canada
        \date{}
            
    \end{center}
\end{titlepage}

\tableofcontents

\pagebreak

\section{Introduction}

The growth in population has resulted in a considerable increase in global energy consumption. This is due to the fact that energy is required for practically all activities like transportation, heating etc. Nowadays, hydrogen has emerged as an attractive energy source that may be utilized to store energy from renewable sources. The time has come to capitalize on hydrogen's potential to play a significant role in addressing critical energy challenges. Recent breakthroughs in renewable energy technology have proved that technical innovation has the potential to build worldwide clean energy firms that do not rely on hydrocarbon-based fuels, therefore helping to reduce pollution levels \cite{3}.

\justify
Because it possesses one of the highest energy densities per unit mass, hydrogen is appealing as a fuel. According to \cite{5}, energy density is defined as the quantity of energy stored per unit volume in a specific system or region of space. According to this statistic, hydrogen has an energy density of 35,000 watts per kilogram (W/kg), but lithium-ion batteries have an energy density of just 200 watts per kilogram (W/kg), and because hydrogen has a tenfold higher energy-to-weight ratio than lithium-ion batteries, it can deliver a greater range while being substantially lighter. Fuel cells, when pure hydrogen is used, can be turned into electricity at high efficiency, durability, and efficiency as needed.

\justify
Because of the physical features of hydrogen, it is difficult to transport and store, prompting us to study safe ways to produce hydrogen locally effectively. One conceivable method for safe hydrogen transportation that appears to be the most efficient is to store the fuel as a low-pressure liquid and then use an onboard ethanol steam reformer (ESR) to produce hydrogen as needed \cite{1}. This is owing to its mobility, renewable nature, and low toxicity.

\justify
The development of real-time efficient and trustworthy control systems to assure device efficiency while reducing the impact of interruptions such as significant variations in interior and external temperatures during transportation, as mentioned in \cite{1}, are crucial for the design of ethanol steam reformers. There have been few mathematically modelled studies on the best design of control techniques for ethanol steam reformers. Previous research has examined the steady-state behaviour of ethanol steam reformers in order to build proportional-integral-derivative (PID)-based decoupled control loops while disregarding physical and operational restrictions and needs \cite{6} and \cite{7}. Another study \cite{4} presented the use of a model which manages a mass flow control of ethanol/water and temperature regulation of a 1 kWe thermal plasma reformer. Although these work, they have certain limitations, such as neglecting physical and operational constraints and not using the system's accessible information. Some of the research relies on linear process models, which might be inaccurate in general, as in \cite{4}. Given the complexity of nonlinear models for ethanol steam reformers, research has been done to employ nonlinear process models, which have a relatively high computational cost for computing the control rule.

\justify
Model predictive control (MPC) is a sophisticated process control method that has been utilized extensively in industrial and chemical processes since the 1980s. The MPC approach is a collection of control techniques that employ a mathematical model of an investigated system to get control actions by minimizing a cost function connected to specified control objectives while taking desired system performance into account \cite{2}.

\justify
The authors of \cite{2} conducted the first research, which was a rigorous examination of the extent of nonlinear model predictive control, in which they analyze a mechanistic model that is a single distributed parameter system which, in fact, makes this physics-based distributed parameter system for ethanol steam reformer a \textit{unique} system. In addition, the same paper presented a method for constructing an approximation reduced-order nonlinear dynamics model with lower computing cost while preserving the same structure of the original equations and physical model parameters.

\justify
Further extending the study \cite{1} introduced a \textit{zero error} approach in the model reformulation. Zero error means that no error is introduced when reformulating the model, and no additional assumptions are made. When putting the model into a nonlinear model predictive control algorithm, this new formulation preserves the same advantages of being dependent on physical model parameters and considerably lowering the real-time computations required. Figure \ref{fig:12} shows the ESR, which is a nonlinear dynamical system comprised of a catalytic ethanol steam reactor in series with a separation stage that includes a selective membrane for hydrogen removal as described in \cite{1}. The mechanistic model described is a function of time and axial direction only (only one spatial dimension) for the system of partial differential equations.

\justify
Our report will provide the first in-depth examination of the aforementioned reformed system, with the goal of attempting to turn it into some type of conservation law or similar standard method ( examples include convection form, diffusion form etc.). We will prove the uniqueness and existence of our system's solutions by understanding their characteristic system, which is usually done by the \textit{Riemann} method. Expanding on this, we look at any potential singularities and how to mitigate this during numerical analysis. This is particularly crucial since any numerical approach we wish to use requires us to discretize the domain (i.e. discretize both the \textit{spatial} and \textit{temporal} parts of the domain and consider a bounded domain). This is significant because it allows us to avoid using numerical approaches that might cause the system to blow up and instead utilize more efficient numerical methods created expressly for 1D conservation rules (examples include The finite volume method, Approximate Riemann Solvers etc).

\justify
Our main contribution with this paper is to perform a numerical comparison of the methods applied to our system based on dependability and efficiency, as defined by computation time and the number of analyses required to achieve a specific level of accuracy, as well as to guarantee that the algorithm is stable, converges and is consistent \cite{3}. This would be accomplished by converting the systems of partial differential equations to a system of ordinary differential equations, for which there are already various techniques for solving this system of ODEs (examples include Euler Methods, Adaptive methods etc). This may also be compared to a direct approach to solving PDE problems utilizing numerical techniques such as finite difference methods, finite element methods, finite volume methods, spectral methods etc. Finally, we take this further by determining the best system formulation and refining the accompanying numerical approach even more to reduce computing costs.

\justify
Aligned with \cite{3}, the goals of the study are to enable a mechanistic model to be employed in real-time control calculations while explicitly accounting for input, state and output constraints with minimal computation cost. This would open up a new field of research into nonlinear distributed parameter systems with additional features such as mass, and particle number, all of which are common in other areas of research (examples include thermodynamics, fluid dynamics etc ). This research may be extended to other reactor systems that are frequently encountered in chemical process control applications etc. Lastly, this would bring us closer to our goal of manufacturing hydrogen safely, which in turn could be used as green energy.

\justify
The remainder of this article is organized as follows: Section 3 describes the ethanol steam reformer model under consideration, as well as the chemical processes and assumptions. Section 4 discusses the original system of partial differential equations as well as the reformed system.
Section 5 defines the conservation law and the convection form of the system. Section 6 summarizes the efforts made to convert the system into conservation legislation. Section 6 illustrates how to solve a first-order quasi-linear PDE problem. The Uniqueness and Existence Theorem for First-Order Linear Systems is covered in Section 7. Sections 8 and 9 explore the Uniqueness and Existence theorem for First Order Quasi Linear Systems and how it cannot be directly applied to our model. Section 10 provides a brief explanation of the numerical methods that are applied to our model. Section 11 summarizes the results as well as the findings and programming methods employed, whereas Section 12 highlights the limits of some of the work done in this research. Finally, Sections 13 and 14 discuss future work and the main conclusions, followed by the bibliography and appendix.

\section{Model Description}

We now proceed to describe the Ethanol Steam Reformer(ESR) as a nonlinear dynamical system that combines in series two process unit operations: 
\begin{enumerate}
	\item a reformer stage comprised of a catalytic ethanol steam reactor
	\item a separation stage comprised of a selective membrane through which hydrogen can penetrate.
\end{enumerate}

\begin{figure}
	\centering
	\includegraphics[width=0.65\linewidth]{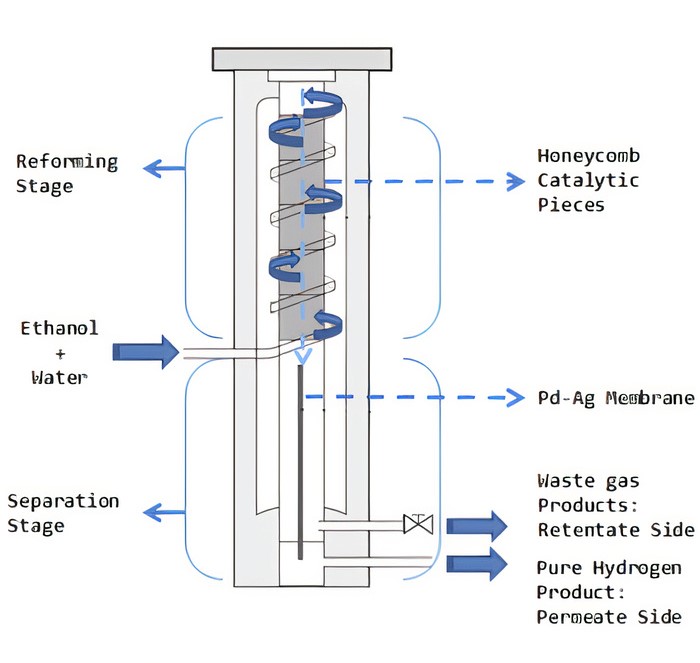}
	\caption{ Staged membrane reactor scheme }
	\label{fig:12}
\end{figure}

\justify
In the reactor, ethanol is reformed with water to generate a gas mixture from which hydrogen is extracted. The entire process takes place within a single integrated module known as a staged membrane reactor which is described in \ref{fig:12}. 

\subsection{Chemical Reactions}
Chemical reactions take place in a tubular packed-bed reactor with a single intake and output. There are 4 primary chemical reactions over cobalt-based catalysts that take place in the staged membrane reactor are as follows \cite{1}

$$\begin{aligned}
	\mathrm{C}_{2} \mathrm{H}_{5} \mathrm{OH} & \longrightarrow \mathrm{CH}_{3} \mathrm{CHO}+\mathrm{H}_{2}, \hspace{20mm} \text{(1a)}\\
	\mathrm{C}_{2} \mathrm{H}_{5} \mathrm{OH} & \longrightarrow \mathrm{CO}+\mathrm{CH}_{4}+\mathrm{H}_{2}, \hspace{18.5mm} \text{(1b)}\\
	\mathrm{CO}+\mathrm{H}_{2} \mathrm{O} & \rightleftharpoons \mathrm{CO}_{2}+\mathrm{H}_{2}, \hspace{30mm} \text{(1c)}\\
	\mathrm{CH}_{3} \mathrm{CHO}+3 \mathrm{H}_{2} \mathrm{O} & \longrightarrow 2 \mathrm{CO}_{2}+5 \mathrm{H}_{2}. \hspace{25mm} \text{(1d)}
\end{aligned}$$
These four reactions are occurring in the same location and under the same circumstances at the same time. To begin, ethanol dehydrogenates to produce hydrogen and formaldehyde (1a), which is then reformed with water to produce carbon dioxide (1d). Furthermore, at normal working circumstances, cobalt catalysts are active for the Water Gas Shift (WGS) process (1c). The unfavourable process is ethanol breakdown, which produces carbon monoxide and methane (1b). The Pd-Ag membrane absorbs only hydrogen during the membrane separation step, leaving waste gases on the retentate side \cite{2}.\\

\justify
For the mathematical modelling of the ethanol steam reformer, various plausible and slightly overlapping assumptions are made which are as follows \cite{1}

\begin{enumerate}
	\item The concentrations in the ethanol steam reformer are easily characterized as two plug-flow units connected in sequence.
	\item The radial dependency of pressure and temperature in the tubular reactor and membrane separator is insignificant.
	\item In the tubular reactor and membrane separator, the fluid is entirely radially mixed at each axial location.
	\item In comparison to convection in the axial direction, the impact of molecular diffusion in the axial direction is modest.
	\item In the operating pressure range, the ideal gas assumption is true.
\end{enumerate}

\justify
These assumptions ensure that the mechanistic model for the two process units is a function of time and axial direction, implying that the partial differential equations describing the system's nonlinear dynamics have just one spatial dimension. This leads to the model being a system of first-order quasi-linear hyperbolic equations.

\section{First Order Quasi Linear PDE System}

We now formulate the above system as a first-order quasi-linear hyperbolic equation. Recall that the class of nonlinear distributed parameter systems described by the first-order quasi-linear hyperbolic vector equation is described as 
$$
A(x, y, u) u_{x}+B(x, y, u) u_{y}=g(x, y, u)
$$
where $u(x, y) \in \mathbb{R}^{m}$ is a vector of distributed states with each element being a function of the real independent variables $x$ and $y, A(x, y, u) \in \mathbb{R}^{m \times m}$ and $B(x, y, u) \in$ $\mathbb{R}^{m \times m}$ are real matrices, $g(x, y, u) \in \mathbb{R}^{m}$ is a vector with each element being an algebraic function of its arguments, and $u_{x}:=\frac{\partial u}{\partial x}$ and $u_{y}:=\frac{\partial u}{\partial y}$. Each individual equation arises from the conservation of some property, such as mass, moles, particles, energy, or momentum. The model equations in \cite{2} can be written as \cite{1}
$$
M f_{t}+N f_{z}=g(t, z, f),
$$
with 
$$
\begin{gathered}
	f:=\left(\begin{array}{c}
		F \\
		T
	\end{array}\right), \quad F=\left(\begin{array}{c}
		F_{1} \\
		\vdots \\
		F_{7}
	\end{array}\right) \\
	M=\left(\begin{array}{cc}
		I-\frac{1}{\|F\|} F e^{\top} & -\frac{1}{T} F \\
		0^{\top} & \frac{c_{v}}{R T\|F\|} F
	\end{array}\right) \\ \\
	N=\left(\begin{array}{cc}
		\frac{R T}{A p} I & 0 \\
		\frac{R T}{A p} e^{\top} & \frac{c_{p}}{A p} F
	\end{array}\right) \\
\\	g =\left(\begin{array}{cc}
	R T [F_1,F_1, \ldots F_1] F V_1R' \\
	R T  [F_2,F_2, \ldots F_2] F V_2R' \\
	\vdots \\
	R T [F_7,F_7, \ldots F_7]F V_7R' \\
	U \beta(T_f - T) - \sum_i F_i [R'H' - V_iR']^{T}
	
\end{array}\right)
\end{gathered}
$$
where
\begin{itemize}

\item $I$ is the $7 \times 7$ identity matrix,
\item $e \in \mathbb{R}^{7}$ is the vector of ones,
 \item $F_{j}$ is the flow rate of species $j$,
\item $\|F\|$ is the Euclidean norm of the vector $F$,
\item $T$ is the temperature (in $\mathrm{K}$ ),
\item $R$ is the ideal gas constant (in $\mathrm{Pa} \mathrm{m}^{3} /(\mathrm{mol} \mathrm{K})$ ),
\item $\rho$ is the membrane thickness (in $\mathrm{m}$ ),
 \item $c_{p}$ and $c_{v}$ are the heat capacities (in $\left.\mathrm{J} /(\mathrm{mol} \mathrm{K})\right)$.
\item $p$ is the pressure (in bar),
\item $A$ is the cross-sectional area of the tubular reactor (in $\left.\mathrm{m}^{2}\right)$,
\item $U$ is the overall heat transfer coefficient (in $\mathrm{J} /\left(\mathrm{m}^{2} \mathrm{~s}\right.$ $\mathrm{K})$ ).
\item $R'$ is a vector of the reactions indices ie $R' = [r_1,r_2,r_3,r_4]$ where $i \in\{1, 2, 3, 4\}$ is the respective reaction index according to the chemical reactions described earlier. 
\item The expressions for the chemical reaction rates are described by elementary kinetics in \cite{2} as follows
$$
\begin{aligned}
	r_{a} &=k_{a} P_{\mathrm{C}_{2} \mathrm{H}_{5} \mathrm{OH}}, \\
	r_{b} &=k_{b} P_{\mathrm{C}_{2} \mathrm{H}_{5} \mathrm{OH}}, \\
	r_{c} &=k_{c}\left(P_{\mathrm{CO}} P_{\mathrm{H}_{2} \mathrm{O}}-\frac{P_{\mathrm{CO}_{2}} P_{\mathrm{H}_{2}}}{k_{W G S}}\right), \\
	r_{d} &=k_{d} P_{\mathrm{CH}_{3} \mathrm{CHO}} P_{\mathrm{H}_{2} \mathrm{O}}^{3}, \\
	k_{W G S} &=e^{\frac{4577.8}{T}-4.33}, \\
	k_{i} &\left.=k_{\infty i} e^{-E_{a_{i}}\left(\frac{1}{R T}-\frac{1}{R T_{\text {rel }}}\right.}\right),
\end{aligned}
$$ 
and where
\item $P_{j}$ is the partial pressure (in bar) of the $j$ th component,
\item $T_{\text {ref }}$ is the reference temperature, equal to $773 \mathrm{~K}$
\item $\beta$ is the area per reactor volume where heat is being transferred (in $\mathrm{m}^{2} / \mathrm{m}^{3}$ ), which is equivalent to $4 / \gamma$, with $\gamma$ being the catalytic reactor diameter (in $\mathrm{m}$ ),
\item $\Delta H_{i}$ is the heat of the $i$ th reaction (in $\mathrm{J} / \mathrm{mol}$ ).
\item $V_i$ is a vector consisting of the stoichiometric coefficient (dimensionless) of all the species in the ith reaction ie $V_i = [\nu_{i, 1},\nu_{i, 2}, \nu_{i, 3}, \nu_{i, 4}]$
\item $H'$ is the vector consisting of $\Delta H_{i}$ ie $H' = [\Delta H_{1}, \Delta H_{2}, \Delta H_{3}, \Delta H_{4}]$
\end{itemize}

\justify
In the same paper, the authors \cite{1} also simplified the above system by writing the vector $f$ in terms of an orthonormal basis with $\|F\|=\sqrt{7} u_{7}$. Define
$$
U:=\left(\begin{array}{c}
	u_{1} \\
	\vdots \\
	u_{7}
\end{array}\right)
$$
and $C_{v}^{u}$ and $C_{p}^{u}$ to be such that
$$
C_{v}^{u}\left(\begin{array}{c}
	u_{1} \\
	\vdots \\
	u_{7}
\end{array}\right)=C_{v} F \text { and } C_{p}^{u}\left(\begin{array}{c}
	u_{1} \\
	\vdots \\
	u_{7}
\end{array}\right)=C_{p} F
$$
Then, 
$$
M e_{i}= \begin{cases}e_{i}, & i=1, \ldots, 6 \\ -\frac{u_{1}}{u_{7}} e_{1}-\cdots-\frac{u_{6}}{u_{7}} e_{6}, & i=7 \\ -\frac{u_{1}}{T} e_{1}-\cdots-\frac{u_{7}}{T} e_{7}+\frac{C_{v}^{u} U}{\sqrt{7} R T u_{7}} e_{8}, & i=8\end{cases}
$$
and
$$
N e_{i}= \begin{cases}\frac{R T}{A p} e_{i}, & i=1, \ldots, 6 \\ \frac{R T}{A p} e_{7}+\frac{\sqrt{7} R T}{A p} e_{8}, & i=7 \\ \frac{C_{p}^{u} U}{A p} e_{8}, & i=8\end{cases}
$$
Hence, the left-hand side of for
$$
u:=\left(\begin{array}{c}
	u_{1} \\
	\vdots \\
	u_{7} \\
	T
\end{array}\right)
$$
can be written as
$$
A u_{t}+B u_{x}  = g(t,x,u)
$$
where
$$
A=\left(\begin{array}{ccc}
	1 & -\frac{u_{1}}{u_{7}} & -\frac{u_{1}}{T} \\
	\ddots & \vdots & \vdots \\
	 1 & -\frac{u_{6}}{u_{7}} & -\frac{u_{6}}{T} \\
	& 0 & -\frac{u_{7}}{T} \\
	& 0 & \frac{C_{v}^{u} U}{\sqrt{7} R T u_{7}}
\end{array}\right)
$$~\\
$$
B=\left(\begin{array}{ccc}
	\frac{RT}{Ap} & 0 & 0  \\
	\ddots & \vdots & \vdots \\
	0 & \frac{RT}{Ap} & 0 \\
0 & \frac{\sqrt 7RT}{Ap} & \frac{C_{v}^{u} U}{Ap}
\end{array}\right)
$$
and
$$g =\left(\begin{array}{cc}
	R T[u_1,u_1, \ldots u_1] U V_1R' \\
	R T [u_2,u_2, \ldots u_2] U V_2R' \\
	\vdots \\
	R T [u_7,u_7, \ldots u_7] U V_7R' \\
	U \beta(T_f - T) - \sum_i u_i [R'H' - V_iR']^{T}
\end{array}\right)
$$\\
The latter model is much easier to work on for numerical analysis due to its simpler matrix structure. This enables us to determine the generalized eigenvalues and left-generalized eigenvectors due to $A$ and $B$ being almost diagonal matrices.
\section{Conservation Law}

For a given physical domain, conservation laws are frequently expressed in integral form \cite{9}. Assume we have a physical domain, $\Omega$, with a domain boundary, $\partial \Omega$. Then, assuming that the physical domain is fixed, the canonical conservation equation is of the type
$$
\frac{d}{d t} \int_{\Omega} U d \mathbf{x}+\int_{\partial \Omega} \mathbf{F}(U) \cdot \mathbf{n} d s=\int_{\Omega} S(U, t) d \mathbf{x}
$$
where $U$ is the conserved state, $\mathbf{F}$ is the flux of the conserved state, $\mathbf{n}$ is the outward pointing unit normal on the boundary of the domain, and $S$ is a source term. Using Gauss's theorem, this conservation law may be represented as a partial differential equation(usually a hyperbolic system) as follows
$$
\int_{\partial \Omega} \mathbf{F} \cdot \mathbf{n} d s=\int_{\Omega} \nabla \cdot \mathbf{F} d \mathbf{x}
$$
$$
\begin{aligned}
	\frac{d}{d t} \int_{\Omega} U d \mathbf{x}+\int_{\partial \Omega} \mathbf{F} \cdot \mathbf{n} d s &=\int_{\Omega} S d \mathbf{x}, \\
	\int_{\Omega} \frac{d U}{d t} d \mathbf{x}+\int_{\Omega} \nabla \cdot \mathbf{F} \cdot \mathbf{n} d s &=\int_{\Omega} S d \mathbf{x}, \\
	\int_{\Omega}\left(\frac{\partial U}{\partial t}+\nabla \cdot \mathbf{F}-S\right) d \mathbf{x} &=0
\end{aligned}
$$
Since this last equation must be valid for any arbitrary domain $\Omega$ and $U_t + \nabla \cdot \mathbf{F}$ is assumed to be regular enough this implies that the integrand must be zero everywhere or equivalently,
$$
\frac{\partial U}{\partial t}+\nabla \cdot \mathbf{F}=S
$$
The above equation is the conservation law written as a partial differential equation.
There are many reasons why converting the system into a conservation law would be useful \cite{9}

\begin{enumerate}
	\item The conservation law's integral form allows for discontinuous solutions (leading to Rankine-Hugoniot relation and correct shock speed). In the non-conservative form, such a character does not exist (it allows only smooth differentiable solutions).
	
	\item Discrete conservation is critical while computing shocks. If the conservation principle is violated, a shock may travel at an incorrect speed.Successful non-conservative procedures may fulfill some form of discrete conservation.
	
	\item Due to the fact that the integral equation might have non-differentiable solutions, the equality of both formulations can break down in some situations, resulting in weak solutions and serious numerical challenges in simulations of such equations.
	
\end{enumerate}

\subsubsection{Example}
\textbf{Burger's Equation}~\\ The inviscid Burgers' equation is a conservation equation, or more precisely, a first order quasilinear hyperbolic equation.
Let $F(u)=\frac{u^2}{2}$ then we have
$$
u_{t}+u u_{x}=0 .
$$

\begin{figure}
	\centering
	\includegraphics[width=0.7\linewidth]{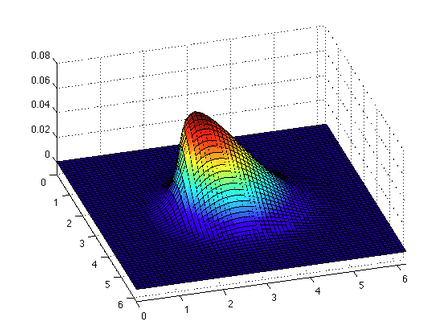}
	\caption{Burgers Equation in two space variables}
	\label{fig:3}
\end{figure}

\subsubsection{Convection}
In Convection we have the almost conservation law form following the fact that if we let $U$ be the 'conserved' scalar quantity and let the fluxes be given by \cite{9},
$$
F=\mathbf{v} U, \quad S=0,
$$
where $\mathbf{v}(\mathbf{x}, t)$ is the known vector of velocity. It should be noted that a non-zero source term might be added but we will proceed with just $S$ being zero. The first order partial differential equation may be used to express this scalar conservation law.
$$
\frac{\partial U}{\partial t}+\nabla \cdot(\mathbf{v} U)=0
$$
$$
\frac{\partial U}{\partial t}+\mathbf{v} \cdot \nabla U+(\nabla \cdot \mathbf{v}) U=0
$$
Now if the velocity field is divergence free $(\nabla \cdot \mathbf{v}=0)$ then we arrive at what is commonly referred to as the convection equation,
$$
\frac{\partial U}{\partial t}+\mathbf{v} \cdot \nabla U=0
$$
Physically, this equation indicates that the quantity $U$ does not change when moving along the stream-wise direction (i.e. convecting with the velocity).

\subsubsection{Transformation}

We take our original model and notice that it breaks up into two sets of equations
\begin{enumerate}
	\item The Molar Conservation Equations
	\item The Energy Conservation Equation
\end{enumerate}

Now the Molar Conservation Equations can easily be converted into a conservation law. For the reforming stage, the molar conservation equations for the seven $j \in(1, \ldots, 7)$ species $\left(\mathrm{C}_{2} \mathrm{H}_{5} \mathrm{OH}, \mathrm{H}_{2} \mathrm{O}, \mathrm{CH}_{4}, \mathrm{H}_{2}\right.$, $\left.\mathrm{CO}, \mathrm{CO}_{2}, \mathrm{CH}_{3} \mathrm{CHO}\right)$ are described by a set of seven nonlinear differential equations, initial conditions, and boundary conditions \cite{1}:
$$
\begin{aligned}
	p \frac{\partial F_{j}}{\partial t}-p & \frac{F_{j}}{\sum_{i} F_{i}} \sum_{i} \frac{\partial F_{i}}{\partial t}-\frac{p}{T} F_{j} \frac{\partial T}{\partial t} = R T \sum_{i} F_{i}\left(\sum_{i} \nu_{j, i} r_{i}\left(\left\{F_{j}\right\}\right)-\frac{1}{A} \frac{\partial F_{j}}{\partial z}\right) \\
	F_{j}(z, 0) &=F_{j, 0}(z), \quad \forall z \in\left[0, \ell_{1}\right] \\
	F_{j}(0, t) &=F_{j, \text { in }}(t), \quad \forall t \geq 0
\end{aligned}
$$
where most of the constants were defined as before in section 4 and the following boundary and initial conditions variables are as follows
\begin{itemize}
	\item - $F_{j}$ is the molar flow (in $\mathrm{mol} / \mathrm{min}$ ) of the $j$ th species,
	\item $F_{j, 0}$ is the molar flow of the $j$ th species at time $t=0$,
	\item $F_{j, \text { in }}$ is the molar flow of the $j$ th species at the reactor inlet $(z=0)$,
	\item $\ell_{1}$ is the axial length (in m) of the reactor,
\end{itemize}
At any given point, the concentration of the $j$ th species, $C_{j}$, and the flow's velocity, $v$, can be computed as \cite{1}
$$
\begin{aligned}
	C_{j} &=\frac{F_{j}}{\sum_{i} F_{i}} \frac{p}{R T} \\
	v &=\frac{1}{A} \sum_{i} F_{i} \frac{R T}{p}
\end{aligned}
$$
Now notice we can transform the Molar Conservation Equation into a conservation equation by using $C_j$ and $v$ as defined above.

$$\begin{aligned}
	&\frac{\partial \mathrm{C}_{\mathrm{j}}}{\partial t}+ \frac{\partial\left(vC_j\right)}{\partial z}=\sum_{i} v_{j, i} r_{i}, \\
	&\mathrm{C}_{\mathrm{j}}(0, z)=\mathrm{C}_{\mathrm{j}, 0}(z), \quad \forall z \in[0,l_1] \\
	&\mathrm{C}_{\mathrm{j}}(t, 0)=\mathrm{C}_{\mathrm{j}, \text { in }}(t), \quad \forall t\geq 0
\end{aligned}$$
which is in a conservation form. The spatiotemporal dynamics of the ESR temperature are described by the energy conservation equation \cite{2}
$$
\begin{array}{r}
	\frac{p}{R T} \frac{\sum_{j} c_{v_{j}} F_{j}}{\sum_{i} F_{i}} \frac{\partial T}{\partial t}=U \beta\left(T_{f}-T\right)-\frac{R T}{A} \sum_{j} \frac{\partial F_{j}}{\partial z}
	-\frac{1}{A}\left(\sum_{j} c_{p_{j}} F_{j}\right) \frac{\partial T}{\partial z}-\sum_{j} \sum_{i} r_{i}\left(\left\{F_{j}\right\}\right)\left(\Delta H_{i}-\nu_{j, i} R T\right),
\end{array}
$$
with initial and boundary conditions
$$
\begin{array}{ll}
	T(z, 0)=T_{0}(z), & z \in\left[0, \ell_{1}\right] \\
	T(0, t)=T_{\text {in }}(t), & \forall t \geq 0
\end{array}
$$
Now the above set of equations was particularly derived from these equations by using $C_j$ and $v$ ie

$$\begin{aligned}
	&\rho_{g} c_{p_{g}}\left(v \frac{\partial T}{\partial z}+\left(1+\frac{\rho_{s} c_{p_{s}}}{\rho_{g} c_{p_{g}}}\right) \frac{\partial T}{\partial t}\right)=U a\left(T_{f}-T\right)+H_{r}\\
	&\begin{aligned}
		&H_{r}=\sum_{j=1}^{N}-\Delta H_{j} r_{j} \quad j=1,2, \ldots, 4 \\
		&T(0, z)=T_{0}(z) \quad z \in[0, L] \\
		&T(t, 0)=T_{i n}(t) \quad \forall t>0
	\end{aligned}
\end{aligned}$$
Now notice that this is not in a conservation form. An attempt was made to convert the above equation into a conservation form. \\

\subsubsection{Convection Form}

Let us now take the above equation and try to convert it into a conservation form as a function of $F^{-1}$ where $F = \sum_{i} F_{i}$ and $v = \frac{FRT}{Ap}$. Let $k = \frac{R}{Ap}$ where $k$ is just a constant then we get $T = \frac{v}{Fk}$. Now substituting this in our equation we get 
\[\rho_{g} c_{p_{g}}\left(v \frac{\partial \frac{v}{Fk}}{\partial z}+\left(1+\frac{\rho_{s} c_{p_{s}}}{\rho_{g} c_{p_{g}}}\right) \frac{\partial \frac{v}{Fk}}{\partial t}\right)=U a\left(T_{f}-T\right)+H_{r} \]

\[= \rho_{g} c_{p_{g}}\left(vk \frac{\partial \frac{v}{F}}{\partial z}+\left(k+\frac{k\rho_{s} c_{p_{s}}}{\rho_{g} c_{p_{g}}}\right) \frac{\partial \frac{v}{F}}{\partial t}\right)=U a\left(T_{f}-T\right)+H_{r}\]

\[ = \rho_{g} c_{p_{g}}kv \frac{\partial \frac{v}{F}}{\partial z} + \rho_{g} c_{p_{g}}k\frac{\partial \frac{v}{F}}{\partial t} +k\rho_{s} c_{p_{s}}\frac{\partial \frac{v}{F}}{\partial t} =U a\left(T_{f}-T\right)+H_{r}\]

\[ = v\frac{\partial \frac{v}{F}}{\partial z} + \frac{\partial \frac{v}{F}}{\partial t} +\frac{\rho_{s} c_{p_{s}}}{\rho_{g} c_{p_{g}}}\frac{\partial \frac{v}{F}}{\partial t} = \left(U a\left(T_{f}-T\right)+H_{r} \right)\frac{1}{\rho_{g} c_{p_{g}}}\frac{1}{k} \]

\[ =  v\frac{\partial \frac{v}{F}}{\partial z} + \left(1 + \frac{\rho_{s} c_{p_{s}}}{\rho_{g} c_{p_{g}}}\frac{\partial \frac{v}{F}}{\partial t} \right) = \left(U a\left(T_{f}-T\right)+H_{r} \right)\frac{1}{\rho_{g} c_{p_{g}}}\frac{1}{k} \]

\justify
Let $\beta = \left(1 + \frac{\rho_{s} c_{p_{s}}}{\rho_{g} c_{p_{g}}}\right)$. Then we get the final equation

\[ =  v\frac{\partial \frac{v}{F}}{\partial z} + \beta\frac{\partial \frac{v}{F}}{\partial t} = \left(U a\left(T_{f}-T\right)+H_{r} \right)\frac{1}{\rho_{g} c_{p_{g}}}\frac{1}{k} \]

\[ =  v\frac{\partial v{F^{-1}}}{\partial z} + \beta\frac{\partial v{F^{-1}}}{\partial t} = \left(U a\left(T_{f}-T\right)+H_{r} \right)\frac{1}{\rho_{g} c_{p_{g}}}\frac{1}{k} \]

\justify
Let $g = v{F^{-1}}$ Then we get 

\[ = v\beta^{-1} \frac{\partial g}{\partial z} + \frac{\partial g}{\partial t} = \left(U a\left(T_{f}-T\right)+H_{r} \right)\frac{1}{\rho_{g} c_{p_{g}}}\frac{\beta^{-1}}{k} \]
\justify
Which is not exactly in the conservation form that we needed but more of in a \textit{convection form}, which is in terms of a function $F^{-1}$.

\subsubsection{Final System}

Finally, collecting both the attempts from above, we have two sets of systems of equations, one of them being

$$\begin{aligned}
	&\frac{\partial \mathrm{C}_{\mathrm{j}}}{\partial t}+ \frac{\partial\left(vC_j\right)}{\partial z}=\sum_{i} v_{j, i} r_{i}, \\
	&\mathrm{C}_{\mathrm{j}}(0, z)=\mathrm{C}_{\mathrm{j}, 0}(z), \quad \forall z \in[0,l_1] \\
	&\mathrm{C}_{\mathrm{j}}(t, 0)=\mathrm{C}_{\mathrm{j}, \text { in }}(t), \quad \forall t\geq 0
\end{aligned}$$
and the other system where $T_f = vF^{-1}_f = g_f$ is just a constant

$$\begin{aligned}
&  v\left(1 + \frac{\rho_{s} c_{p_{s}}}{\rho_{g} c_{p_{g}}}\right)^{-1}\frac{\partial \left(v{F^{-1}}\right)}{\partial z} + \frac{\partial \left(v{F^{-1}}\right)}{\partial t} = \left(U a\left(F^{-1}_f-\frac{v}{Fk}\right)+H_{r} \right)\frac{1}{\rho_{g} c_{p_{g}}}\frac{Ap}{R} \left(1 + \frac{\rho_{s} c_{p_{s}}}{\rho_{g} c_{p_{g}}}\right)^{-1}\\
&H_{r}=\sum_{j=1}^{N}-\Delta H_{j} r_{j} \quad j=1,2, \ldots, 4 \\
&vF^{-1}(0, z)=vF^{-1}_{0}(z) \quad z \in[0, L] \\
&vF^{-1}(t, 0)=vF^{-1}_{i n}(t) \quad \forall t>0
\end{aligned}$$
which is equivalent to

$$\begin{aligned}
& \frac{\partial g}{\partial t} + v\beta^{-1}\frac{\partial \left(g\right)}{\partial z} = \left(U a\left(g_f-\frac{g}{k}\right)+H_{r} \right)\frac{1}{\rho_{g} c_{p_{g}}}\frac{1}{k} \beta^{-1}\\
&H_{r}=\sum_{j=1}^{N}-\Delta H_{j} r_{j} \quad j=1,2, \ldots, 4 \\
&g(0, z)=g_{0}(z) \quad z \in[0, L] \\
&g(t, 0)=g_{i n}(t) \quad \forall t>0
\end{aligned}$$
As a result, our PDE system would be a tuple $\left(C_1,C_2,C_3,C_4,C_5,C_6,C_7,g\right)$, which is not entirely in the conservation form but rather a hybrid of the conservation and convection forms. As a consequence, the entire system is not in conservation form. We could further decompose it by decomposing the 8 by 8 system into the 7 by 7 system of Conservation form and the last one on its own, but because this is a coupled first-order partial differential equation, we would have to treat the system as a whole rather than independently.

\section{Solution of the System}
We now proceed to derive the solution of the above system. We use the method of characteristics to help aid us in solving this, which is introduced as follows

\subsection{Method of Characteristics}
The primary concept of the method of characteristics is to convert a partial differential equation to an ordinary differential equation using characteristic curves, which enables us to solve the PDE.\\
A parametric curve is a the image of a mapping $\gamma(t) \in \mathbb{R}^{2}$ where $t$ ranges in some open interval, $\gamma(t):=(x(t), y(t))$. Assume $u=u(x, y)$ is an unknown function in the $(x, y)$-plane \cite{13}. The restriction of $u$ to $\gamma$ is just $u(\gamma(t)):=u(x(t), y(t))$. An ordinary differential equation of $u$ along $\gamma$ if of the form
$$
\frac{d u}{d t}(\gamma(t))=f(t, u(\gamma(t)))
$$
where $\frac{d u}{d t}(\gamma(t))$ is the derivative of $u$ along $\gamma(t)$.
The characteristic technique is based on locating a family of curves known as the characteristic curves associated with a particular partial differential equation, such that the partial differential equation reduces to an ordinary differential equation along the curves which then enables us to solve the PDE. 
For example, consider the family of curves $\gamma(t ; c):=(t+c, t)$ for a parameter $c$. If $u=u(x, y)$ is a continuously differentiable function, then $\frac{d}{d t} u(\gamma(t))$ by the chain rule has the following form
$$
\frac{d u}{d t}(\gamma(t))=\partial_{x} u \frac{d x}{d t}+\partial_{y} u \frac{d y}{d t}=\partial_{x} u+\partial_{y} u
$$
Therefore, the partial differential equation $\partial_{x} u+\partial_{y} u=0$ reduces to the ordinary differential equation $\frac{d u}{d t}=0$ along $\gamma(t)$. The solution is
$$
u(\gamma(t))=u(\gamma(0))=u(c, 0),
$$
for the parameter $c$. If $u(c, 0)$ is provided, for example, $u(c, 0)=h(c)$, then by the relation $c=x-y$, we obtain $u(x, y)=h(x-y)$ that is the general solution of the equation. This can easily be extended to systems of PDE by considering the characteristic lines in a vectorized format, which will be presented in the solution of our above system.

\subsection{Solution}

This goes as follows \cite{1}. The first step is to multiply both sides by a row vector $\boldsymbol{l}^{\top}$ with $\boldsymbol{l} \in \mathbb{R}^{m}$ :
$$
\boldsymbol{l}^{\top} u_{t}+\boldsymbol{l}^{\top} B u_{x}=\boldsymbol{l}^{\top} g .
$$
The method of characteristics analyzes the dynamics in terms of a characteristic curve, which is the ordered pair $(x(\omega), t(\omega))$ parameterized by the real scalar $\omega$. Multiplying each term on the left-hand side by $x(\omega)$ and $t(\omega)$ gives
$$
\left(\frac{1}{t_{\omega}} \boldsymbol{l}^{\top}\right) u_{t} t_{\omega}+\left(\frac{1}{x_{\omega}} \boldsymbol{l}^{\top} B\right) u_{x} x_{\omega}=\boldsymbol{l}^{\top} g .
$$
Now if a vector $l$ and ordered pair $(x(\omega), t(\omega))$ can be found such that
\[
\frac{1}{t_{\omega}} \boldsymbol{l}^{\top} =\frac{1}{x_{\omega}} \boldsymbol{l}^{\top} B,
\label{1} \tag{1} \]
then this can be written as
$\left(\frac{1}{t_{\omega}} \boldsymbol{l}^{\top} \right)\left(u_{x} x_{\omega}+u_{t} t_{\omega}\right)=\boldsymbol{l}^{\top} g .$
Application of the chain rule and multiplication by a scalar simplifies $u_\omega = u_x x_\omega + u_t t_\omega$ to
$\left(\boldsymbol{l}^{\top}u_\omega\right)=\boldsymbol{l}^{\top} g t_\omega.$ Then this gives the equivalent expression 
\[\left(\boldsymbol{l}^{\top}Bu_\omega\right)=\boldsymbol{l}^{\top} g x_\omega\]
The above expressions can be written using the chain rule as
\[\quad l^{\top} \frac{d u}{d t}=l^{\top} g\] and
\[\quad l^{\top} B \frac{d u}{d x}=l^{\top} g\]
The next step is to characterize the set of vectors $l$ and ordered pairs $(x(\omega), t(\omega))$ that satisfy the equation. From \ref{1} we define the left generalized eigenvalue decomposition; that is, multiplication by $t_{\omega}$ gives the equivalent expression
$\boldsymbol{l}^{\top} B \frac{d u}{d t}=\boldsymbol{l}^{\top} g$ characterize the $t(\omega))$ that satisfy generalized eigenvalues by $x_{\omega}$ gives the $\boldsymbol{l}^{\top}(B-\sigma I)=0$
where
$$
\sigma=\frac{x_{\omega}}{t_{\omega}}=\frac{d x}{d t}
$$
The following equation has $m$ solutions (possibly repeated) defined by the generalized eigenvalues $\sigma_{k}$ and left generalized eigenvectors $\boldsymbol{l}_{k}$
$$
l_{k}^{\top}\left(B-\sigma_{k} I\right)=0, \quad k=1,2, \ldots, m
$$ The hyperbolicity of $\boldsymbol{l}$ is defined by the requirement that all the roots $\sigma_k, \; \forall k \in 1,2 \ldots$ of $(B-\sigma I)=0$ are real and that in addition there exists $m$ linearly independent eigenvectors. The last step is to collect the above set of equations into a system of ODEs,
$$
\begin{aligned}
	\frac{\mathrm{d} x}{\mathrm{~d} t} &=\frac{1}{\sigma_{k}(x, t, u)} \\
	\boldsymbol{l}_{k}(x, t, u)^{\top} A(x, y, u) \frac{\mathrm{d} x}{\mathrm{~d} t} &=\boldsymbol{l}_{k}(x, t, u)^{\top} g(x, t, u) \\
	k=1,2, \ldots, m.
\end{aligned}
$$
which can easily be solved by standard techniques either by hand or numerically.

\section{Uniqueness and Existence Theorem}
	
\subsection{Our Model}
We have our system of Partial Differential Equations, which is given by 
\[Au_t + Bu_x = g(x,t,u)\] Now a problem arises here as we have the fact that $\operatorname{det}A = 0$ and so we cannot transform it into the equation
\[u_t + Bu_x = g(x,t,u)\]
We will thereby first prove the uniqueness and existence theorem to the above system and then show how it can be applied to our system.

\subsection{Uniqueness Theorem for Linear Systems of First Order}
The method of "energy integrals" permits a transparent uniqueness proof for hyperbolic systems of first order \cite{14}. We write the in- homogeneous system in the vector and matrix form
\[
u_{t}+A u_{x}+B u+c=0 \label{2} \tag{2}
\]
where $u$ is the unknown vector function, $A(x, t)$ and $B(x, t)$ are given matrices, and $c(x, t)$ is a given vector. We assume, moreover, that $A$ has continuous derivatives ie ($A \in C^{1}$ ) and that $B$ and $c$ are continuous. Without the restriction of generality, we assume for the proof that $A=\left(a_{i k}(x, t)\right)$ is a symmetric matrix, and we can convert a linear hyperbolic system into a symmetric form (e.g., into the characteristic normal form). The following is done for a quasi-linear hyperbolic system( Similarly follows for a linear system)

\subsubsection{Symmetric Form}

We have that \[u_{t}+B(x, t, u) u_{x}+g(x, t, u)=0 \] is a first-order quasi linear system. We now set $U=H u$, where the rows of the matrix $H$ are the eigenvectors $l_{k}$, and substitute
$$
u=H^{-1} U, \quad u_{x}=H^{-1} U_{x}+H_{x}^{-1} U, \quad u_{t}=H^{-1} U_{t}+H_{t}^{-1} U
$$
into \eqref{1}. Then we obtain
$$
B H^{-1} U_{x}+H^{-1} U_{t}= -g -\left(B H_{x}^{-1}+H_{t}^{-1}\right) U .
$$
Now, multiplying by $H$ from the left, we obtain the diagonal form
$$
U_{t}+H B H^{-1} U_{x}=G,
$$
where $G=H\left(-g-B H_{x}^{-1} U-H_{t}^{-1} U\right)$ is a known vector, depending linearly on $U$ but not depending on the derivatives of $U$. Thus setting $B' = HBH^{-1}$ yields the final equation 
\[U_t + B'U_x = G\]
which is in the symmetric form.

\subsubsection{Characteristic Lines}

\begin{figure}
	\centering
	\includegraphics[width=0.5\linewidth]{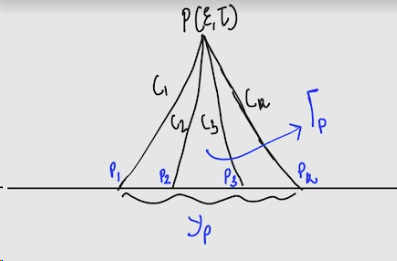}
	\caption{Characteristic Lines}
	\label{fig:1}
\end{figure}

Through a point $P$ with coordinates $\left(\xi, \tau\right)$, we draw the characteristics $C_{1}, C_{2}, \cdots, C_{k}$ backward, obtaining the points $P_{1}, P_{2}, \cdots, P_{k}$ of intersection with the initial line $t=0$. This is drawn in \ref{fig:1}

\subsubsection{Uniqueness Theorem Restated}
\begin{theorem}
	If $c=0$ and $u(x, 0)=0$ in a closed interval of the line $t=0$ including all the points $P_{x}$, then $u(\xi, \tau)=0$.
\end{theorem}
The smallest such interval, i.e., the interval cut out by the outer characteristics through $P$, contains the domain of dependence of $P$. As a matter of fact, if $u(x, 0)=0$ on $[P_{1}, P_{k}]$, then $u(x, t)=0$ in the triangular region $\Gamma_{P}$ \ref{fig:1} cut out by the outer characteristics through $P$ and having the interval $[P_{1}, P_{k}]$ of the line $t=0$ as base. 

\subsubsection{Proof}

We use Green's identity, which is equivalent to
$$
(u, A u)_{x}=\left(u_{x}, A u\right)+\left(u, A_{x} u\right)+\left(u, A u_{x}\right)
$$
since $A$ is symmetric, $\left(u, A u_{x}\right)=\left(A u, u_{x}\right)$, this identity reduces to
$$
2\left(u, A u_{x}\right)=(u, A u)_{x}-\left(u, A_{x} u\right) .
$$
Taking the inner product of the differential equation \eqref{2} with the vector $u$ and using the preceding relation we have (for $c=0$ )
$$
\frac{1}{2}(u, u)_{t}+\frac{1}{2}(u, A u)_{x}-\frac{1}{2}\left(u, A_{x} u\right)+(u, B u)=0 .
$$
We now introduce with a constant $\mu$ instead of $u$ the unknown vector
$$
v=e^{-\mu t} u
$$
which leads to the differential equation
\[\quad v_{t}+A v_{x}+ \left(B +\mu I\right) v=0 \tag{3} \label{3}\]
and the initial condition $v(x, 0)=0$. Let us denote $B*$ with the unit matrix $I$ to be the following,
$$
B^{*}=B+\mu I .
$$
The preceding quadratic identity becomes if we write $u$ again instead of $v$,
$$
\frac{1}{2}(u, u)_{t}+\frac{1}{2}(u, A u)_{x}=(u, \hat{B} u) .
$$
The quadratic form on the right-hand side is formed with the matrix
$$
\hat{B}=-B^{*}-\frac{1}{2} A_{x},
$$
which is the aim of the construction of our constant $\mu$. By choosing $\mu$ sufficiently large, we can ensure that $(u, \hat{B} u) \leq 0$.\\ 
\begin{figure}
	\centering
	\includegraphics[width=0.5\linewidth]{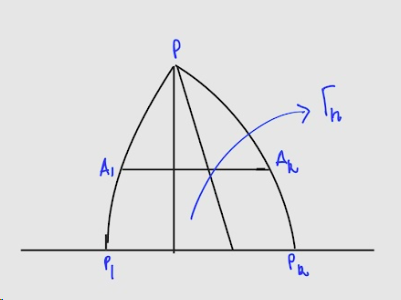}
	\caption{Trapezoid region}
	\label{fig:1234}
\end{figure}

\justify
Now we integrate over the trapezoid-like domain $\Gamma_{h}$ in \ref{fig:1234} which has the boundary $\beta_{h}=P_{1} P_{k} A_{k} A_{1}$ denoting the components of the outward unit normal vector by $x_{\nu}, t_{\nu}$, and obtain Green's formula
$$
\begin{aligned}
	&\frac{1}{2} \iint_{\Gamma_{h}}\left[(u, u)_{t}+(u, A u)_{x}\right] dx dt \leq 0 \\
	&=\frac{1}{2} \int_{\Gamma_{h}}\left[(u, u) t_{\nu}+(u, A u) x_{\nu}\right] ds \leq 0 \\
	&=\frac{1}{2} \int_{A_{1} A_{k}+P_{k} P_{1}}(u, u) d x+\frac{1}{2} \int_{C_{1}+C_{k}} x_{\nu}\left(u,\left[A+\frac{t_{\nu}}{x_{\nu}} I\right] u\right) d s
\end{aligned}
$$
With the notation
$$
E(h)=\frac{1}{2} \int_{A_{1}}^{A_{k}}(u, u) d x
$$
this means
\[E(h)-E(0) \leq-\frac{1}{2} \int_{C_{1}+C_{k}} x_{\nu}\left(u,\left[A+\frac{t_{\nu}}{x_{\nu}} I\right] u\right) ds \tag{4} \label{4} \]

We now show that the right-hand side is non-positive. For this purpose we recall that $C_{1}$ and $C_{k}$ are characteristic curves $\left(\phi^{1}(x, t), \ldots \phi^{k}(x, t)\right) = \vec{0}^{t}$ and that, therefore, they satisfy the differential equations
$$
-[\frac{\phi_{t}^{1}}{\phi_{x}^{1}} , \dots \frac{\phi_{t}^{k}}{\phi_{x}^{k}}] = [l_1, \ldots l_{k}]
$$
where $\left(\l_{1}, \dots \l_{k}\right)$ are eigenvalues of the matrix $A$, i.e., values for which $A-\l I$ is singular.

\justify
Now along $C_{1}=\phi^{1}$, the outward normal has components proportional to $\phi_{x}^{1}, \phi_{t}^{1}$; hence $-t_{\nu} / x_{\nu}=\tau^{1}$. Since $C_{1}$ is the outer characteristic at the left, $\l_{1}$ is the largest eigenvalue of $A$. Similarly, along $C_{k}=\phi^{k}$, the components of the normal are proportional to $\phi_{x}^{k}, \phi_{t}^{k}$ and $-t_{v} / x_{v}=\l_{k}$ is the smallest eigenvalue of $A$. We recall that, for a symmetric matrix $A$, the extreme eigenvalues can be characterized by
$$
\l_{1}=\max \frac{(u, A u)}{(u, u)}, \quad \l_k=\min \frac{(u, A u)}{(u, u)}
$$Consequently
$$
(u, u) \l_{1} \geq(u, A u), \quad(u, u) \l_{k} \leq(u, A u)
$$
or
$$
\left(u,\left[A-\l_{1} I\right] u\right) \leq 0, \quad\left(u,\left[A-\l_k I\right] u\right) \geqq 0 .
$$
Since $x_{\nu}<0$ on $C_{1}$ and $x_{\nu}>0$ on $C_{k}$, we have
$$
-\frac{1}{2} \int_{C_{i}} x,\left(u,\left[A+\frac{t_{\nu}}{x_{\nu}} I\right] u\right) d s \leq 0 \quad \forall i \in \{1,\ldots k\}
$$
This shows that the right-hand side of \ref{4} is non-positive, and thus we can write it as
$$
E(h) \leq E(0) .
$$
Since by assumption $E(0)=0$, we have $0 \leq E(h) \leq 0$; hence $E(h)=0$ for $h>0$; therefore $u=0$ for $t=h$. Hence the uniqueness theorem is proved for First Order Linear Systems.

\section{Uniqueness Theorem for Quasi Linear Systems of First Order}
The following uniqueness theorem remains valid for quasi-linear systems of first order \cite{14}
\[u_{t}+B(x, t, u) u_{x}+g(x, t, u)=0 \label{5} \tag{5}\]
even though the characteristic lines $C_{k}$ of the above equation depend on the solution $u$.

\justify
We assume here that the matrices $B$ and $g$ possess continuous derivatives with respect to $x, t$, and $u$ in the region under consideration. Let $u$ and $v$ be two solutions of the above system defined in a domain $D$ and having the same initial values on the initial interval. Then we consider the following equations
\begin{align}
	& u_{t}+B(x, t, u) u_{x}+g(x, t, u)=0 \\ 
	& u(x,0) = \psi(x)
\end{align} and

\begin{align}
	& v_{t}+B(x, t, v) v_{x}+g(x, t, v)=0 \\ 
	& v(x,0) = \psi(x)
\end{align}
\justify
Subtracting the differential equation for $v$ from that for $u$ and denoting this by $z$ ie ($z(x,t) = u(x,t) - v(x,t)$) we get
\begin{align}
	& z_{t}+B(v) z_{x}+[B(u)-B(v)] u_{x}+g(u)-g(v)=0 \\
	& z(x,t) = 0
\end{align}
Because both $B$ and $g$ are differentiable and continuous, we may apply the mean value theorem
$$
B(u)-B(v)=H(u, v) z ; \quad g(u)-g(v)=K(u, v) z
$$
where $H, K$ are continuous functions. We now consider $u, v, u_{x}$ as known expressions in $x, t$ and substitute these expressions in $H$ and $K$ as well as in $B(v)$; thereby becomes a linear homogeneous differential equation for $z$ of the form
$$
z_{t}+b z_{x}+g z=0
$$
with initial values zero. The theorem proved in the above section asserts now that $z$ is identically zero, and thus the uniqueness theorem is proved also for quasi-linear equations.

\section{Existence Theorem for Quasi Linear Systems of First Order}

We shall now prove the existence theorem for first-order Quasi-linear systems
\[u_{t}+B(x, t, u) u_{x}+g(x, t, u)=0\]
\justify
Note that the characteristics $C_{k}$, eigenvalues $\sigma_k$ and eigenvectors $l_k$ $\forall k \in \{1,2, \dots m\}$  now depend on the specific function $u$ \cite{14}.
\subsubsection{Assumptions}
\begin{enumerate}
	\item The matrix $B(x, t ; v)$ has $k$ real eigenvalues $\sigma_k(x, t ; v)$, briefly $\sigma(v)$, and a matrix $\Lambda(v)$ of independent eigenvectors $l_{k}$ with Lipschitz-continuous derivatives with respect to all arguments in a fixed domain $G_{h}$ to be described presently provided we substitute for $v$ any "admissible" function restricted by the above conditions
	\item The matrix $g$, and the initial values $\psi(x)$ possess Lipschitz-continuous first derivatives with respect to $x, t, u$, then in a suitable neighbourhood $0 \leq t \leq h$ of a portion $g$ of the $x$-axis a uniquely determined solution with Lipschitz-continuous first derivatives exists, provided that the system is hyperbolic for the given initial values $u(x, 0)=\psi(x)$.
	\item We now consider functions $v$, not necessarily solutions, with the fixed initial values $v(x, 0)=$ $\psi(x)$, with Lipschitz-continuous first derivatives and satisfying the inequalities
	$$\|v\| \leq M, \quad \max\{|v_x|,|v_t|\} \leq M_{1}$$
	with fixed $M$ and $M_{1}$. 
\end{enumerate}

\subsubsection{Proof}

The solutions of the ordinary differential equations $d x / d t=\sigma_{k}$ are called the characteristics $C_{k}$ of the $v$-field. For admissible functions, their slopes are uniformly bounded: $\left|\sigma_k\right|<\mu$, and we now specify that a closed domain $G$ and in it the strip $G_{h}$ consists of points such that the $v$-characteristics $C^{*}$. followed backwards towards $t=0$ remain in $G$ and intersect the portion $g$ of the $x$-axis. The set or "space" of all such functions $v$ is again called $S_{h}$. The eigenvectors $l$ and eigenvalues, as well as their derivatives, depend on the variables $x, t$-and $v$-and are Lipschitz continuous.

\justify
We define the characteristic differential operators in the $v$-field by $D=\partial / \partial t+l(\partial / \partial x)$. Let us now proceed by a natural iteration scheme.
\subsubsection{Iteration}

To allow enough leeway for $v$, we may set $M=N+1$ assuming a bound $\|\psi(x)\|<N$. Then, after substituting an admissible function $v(x, t)$ in $B$ and $g$, the above equation becomes linear. The solution $u$ of Cauchy's problem with the given initial values $u(x, 0)=\psi(x)$ is called $u=T v$, and we obtain the solution $u$ of the above equation as a fixed element of the transformation $T$, specifically as the uniform limit for $n \rightarrow \infty$ of iterations $T u_{n}=u_{n+1}$ with $u_{0}=\psi(x)$, where $T=T_{n}$ depends on $n$.

\justify
First, we obtain $u=T V$ according to the procedure, based on the introduction of a new function $U=\Lambda u, V=\Lambda v, \Psi=\Lambda \psi$, etc. Then in the basic formulas in \cite{1}, we have merely to observe that the differentiation with respect to $x$ must account for the dependence of $l, \cdots$, on $v$; that is, we must interpret $(d / d x) l=l_{x}+l_{v} v_{x}, \cdots$, thus introducing derivatives of $v$. This leads quite directly to the lemma: 

\subsubsection{Lemma}For $M=$ $N+1$ we can choose a sufficiently large bound $M_{1}$ and sufficiently small $h$ such that any function $v$ of $S_{h}$ is transformed into another function $w=T v$ of $S_{h}$.

\justify
Furthermore, the sequence $u_{n}$ converges uniformly in $G_{h}$ to a limit function $u$. This function is uniquely determined and obviously satisfies the differential equation in the characteristic form $l D u+$ $l g=0$ where 

\[D =\partial / \partial t+l(\partial / \partial x)\]
\justify
and $l$ can be chosen $l = l_k$ along any characteristic line $C_k$.
The contracting property and uniqueness follow directly by considering the difference $z=u-u^{*}$ of two admissible functions $u=T v$ and $u^{*}=T v^{*}$ with the difference $\zeta=v-v^{*}$, the differential equation
$$
z_{t}+B(v) z_{x}+\left(B(v)-B\left(v^{*}\right)\right) z_{x}^{*}+g(v)-g\left(v^{*}\right)=0 .
$$
Using the mean value theorem, we have
$$
z_{t}+B(v) z_{x}+\zeta K=0,
$$
where $K$ is a bounded function of $x, t$, and obtain an inequality of the form
$$
\|z\|<M_{3} h k\|K\| \text {. }
$$
With sufficiently small $h$ we assure the inequality $\|z\| \leq \frac{1}{2}\|\zeta\|$, and hence the contracting character of the transformation $T$. Therefore $u_{n}$ converges uniformly to a function $u$ in $G_{h}$. Hence the uniqueness theorem is proved.

\subsection{Theoretical Analysis of our System}
We have the reformulated system as \[Au_t + Bu_x = g(x,t,u)\] Now, as stated previously; a problem arises here as we have the fact that $\operatorname{det}A = 0$, and so we cannot transform it into the equation
\[u_t + Bu_x = g(x,t,u)\] 
\justify
An attempt was made to prove the uniqueness theorem for the system by first trying to first collect the set of equations $u = (u_1,u_2,u_3,u_4,u_5,u_6,T)$ into a 7 by 7 system. We then try to rewrite $u_7$ in terms of the other $u_i, \forall i \in {1,\ldots 6}$ and $T$. This enables $u_7$ to be written in terms of an integral form, and then an analysis could be made, but due to lack of time, we could not proceed. Hence we take a different approach and try to use the original system formulation ie
$$\begin{aligned}
	&M f_{t}+N f_{z}=g(t, z, f)\\
	&f:=\left(\begin{array}{c}
		F \\
		T
	\end{array}\right), \quad F=\left(\begin{array}{c}
		F_{1} \\
		\vdots \\
		F_{7}
	\end{array}\right)\\
	&M=\left(\begin{array}{cc}
		I-\frac{1}{\|F\|} F e^{\top} & -\frac{1}{T} F \\
		0^{\top} & \frac{C_{v}}{R T\|F\|} F
	\end{array}\right) \text {, }\\
	&N=\left(\begin{array}{cc}
		\frac{R T}{A p} I & 0 \\
		\frac{R T}{A p} e^{\top} & \frac{C_{p}}{A p} F
	\end{array}\right)
\end{aligned}$$ 

\justify
with initial and boundary conditions as follows

$$\begin{aligned}
	&F_{j}(z, 0)=F_{j, 0}(z), \quad \forall z \in\left[0, \ell_{1}\right]\\
	&F_{j}(0, t)=F_{j, \text { in }}(t), \quad \forall t \geq 0\\
	&T(z, 0)=T_{0}(z), \quad z \in\left[0, \ell_{1}\right], \\
	& T(0, t)=T_{\text {in }}(t), \quad \forall t \geq 0
\end{aligned}$$

\justify
Notice we have that $\operatorname{det}(M) \neq 0$ which was also confirmed by using Mathematica. Therefore we may solve the above system with respect to $u_t$ and write it in the equivalent form

\[f_t + M^{-1}Nf_z = M^{-1}g(z,t,f)\]
\[= f_t + N'f_z = g'(z,t,f)\]
with $N' = M^{-1}N$ and $g' =M^{-1}g(z,t,f)$, which implies the line $t = \text{constant}$ are now non characteristic or free. Hence by using this system, we just apply the above uniqueness and existence theorem that was proved earlier to the above system. 
\section{Numerical Analysis}

We take the reformulated system $Au_t + Bu_x$ to obtain the generalized eigenvalues and left eigenvectors as mentioned in Section 6. Notice over here that we need to find the generalized eigenvalues and left eigenvectors of $\operatorname{det}(A-\sigma B)$. This is solved in \cite{1}, and the solution is as follows.
In particular, the analytical solution to $\operatorname{det}(A-\sigma B)=0$ is
$$
\begin{gathered}
	\sigma_{k}=\frac{A p}{R T}, \quad k=1, \cdots, 6, \\
	l_{1}=\left(\begin{array}{c}
		1 \\
		0 \\
		\vdots \\
		0
	\end{array}\right), \ldots, \boldsymbol{l}_{6}=\left(\begin{array}{c}
		0 \\
		\vdots \\
		0 \\
		1 \\
		0 \\
		0
	\end{array}\right)
\end{gathered}
$$
and $\sigma_{7}$ and $\sigma_{8}$ solve
$$
\operatorname{det}\left(\begin{array}{cc}
	-\sigma \frac{R T}{A p} & -\frac{u_{7}}{T} \\
	-\sigma \frac{\sqrt{7} R T}{A p} & \frac{C_{v}^{u} U}{\sqrt{7} R T u_{7}}-\sigma \frac{C_{p}^{u} U}{A p}
\end{array}\right)=0
$$
which can be written as
$$
\frac{R T C_{p}^{u} U}{(A p)^{2}} \sigma^{2}-\left(\frac{C_{v}^{u} U}{\sqrt{7} A p u_{7}}-\frac{\sqrt{7} R u_{7}}{A p}\right) \sigma=0
$$
whose solutions are
$$
\sigma_{7}=0, \quad \sigma_{8}=\frac{A p}{R T C_{p}^{u} U}\left(\frac{C_{v}^{u} U}{\sqrt{7} u_{7}}-\sqrt{7} R u_{7}\right)
$$
Since the eigenvalues are known, computing the left generalized eigenvectors $l_{7}$ and $l_{8}$ is straightforward. We obtain the following set of eigenvectors for $l_7$ and $l_8$ when computed, this can be done numerically or can be found out as follows

\[l_7^{T}A = 0\]
\[(l_7^{1},l_7^{2},l_7^{3},l_7^{4},l_7^{5},l_7^{6},l_7^{7},l_7^{8})\left(\begin{array}{ccc}
	1 & -\frac{u_{1}}{u_{7}} & -\frac{u_{1}}{T} \\
	\ddots & \vdots & \vdots \\
	1 & -\frac{u_{6}}{u_{7}} & -\frac{u_{6}}{T} \\
	& 0 & -\frac{u_{7}}{T} \\
	& 0 & \frac{C_{v}^{u} U}{\sqrt{7} R T u_{7}}
\end{array}\right) = 0\]
\[(l_7^{1},l_7^{2},l_7^{3},l_7^{4},l_7^{5},l_7^{6},l_7^{7},l_7^{8}) = \left\{\frac{{u}_1}{{u}_7}, \frac{{u}_2}{{u}_7}, \frac{{u} 3}{{u}_7}, \frac{{u}_4}{{u}_7}, \frac{{u}_5}{{u}_7}, \frac{{u}_6}{{u}_7}, 1,0\right\}\] 
which implies
$$\boldsymbol{l}_{7}=\left(\begin{array}{c}
\frac{{u}_1}{{u}_7}\\
 \frac{{u}_2}{{u}_7}\\
  \frac{{u} 3}{{u}_7}\\
   \frac{{u}_4}{{u}_7}\\
    \frac{{u}_5}{{u}_7}\\
     \frac{{u}_6}{{u}_7}\\
      1\\
      0
\end{array}\right)$$
\justify
Similarly, we can find out the left eigenvector $l_8$ as follows and its given by 

\[l_8^{T}(A-\sigma_8B) = 0\]
\[(l_8^{1},l_8^{2},l_8^{3},l_8^{4},l_8^{5},l_8^{6},l_8^{7},l_8^{8})\left(\left(\begin{array}{ccc}
	1 & -\frac{u_{1}}{u_{7}} & -\frac{u_{1}}{T} \\
	\ddots & \vdots & \vdots \\
	1 & -\frac{u_{6}}{u_{7}} & -\frac{u_{6}}{T} \\
	& 0 & -\frac{u_{7}}{T} \\
	& 0 & \frac{C_{v}^{u} U}{\sqrt{7} R T u_{7}}
\end{array}\right) - \frac{A p}{R T C_{p}^{u} U}\left(\frac{C_{v}^{u} U}{\sqrt{7} u_{7}}-\sqrt{7} R u_{7}\right)\left(\begin{array}{ccc}
\frac{RT}{Ap} & 0 & 0  \\
\ddots & \vdots & \vdots \\
0 & \frac{RT}{Ap} & 0 \\
0 & \frac{\sqrt 7RT}{Ap} & \frac{C_{u}^{p} U}{Ap}
\end{array}\right)\right)= 0\]

$$\boldsymbol{l}_{8}=\left(\begin{array}{c}
	\frac{7 u_1\left(-C_{u}^{v} C_{u}^{p} UU^{\top} u_7+\sqrt{7}{C_{u}^{p}}^2 UU^{\top} u_7^{2} - 7 C_{u}^{p} R U u_7^{3}\right)}{T\left(C_{u}^{v} U + 7 R u_7^{2}\right)\left(\sqrt{7} C_{u}^{v}U - 7 C_{u}^{p} U u_7 + 7 \sqrt{7} R u_7^{2}\right)}\\~\\
	
	\frac{7 u_2\left(-C_{u}^{v} C_{u}^{p} UU^{\top} u_7+\sqrt{7}{C_{u}^{p}}^2 UU^{\top} u_7^{2} - 7 C_{u}^{p} R U u_7^{3}\right)}{T\left(C_{u}^{v} U + 7 R u_7^{2}\right)\left(\sqrt{7} C_{u}^{v}U - 7 C_{u}^{p} U u_7 + 7 \sqrt{7} R u_7^{2}\right)}\\~\\
	
	\frac{7 u_3\left(-C_{u}^{v} C_{u}^{p} UU^{\top} u_7+\sqrt{7}{C_{u}^{p}}^2 UU^{\top} u_7^{2} - 7 C_{u}
		^{p} R U u_7^{3}\right)}{T\left(C_{u}^{v} U + 7 R u_7^{2}\right)\left(\sqrt{7} C_{u}^{v}U - 7 C_{u}^{p} U u_7 + 7 \sqrt{7} R u_7^{2}\right)}\\~\\
	\frac{7 u_4\left(-C_{u}^{v} C_{u}^{p} UU^{\top} u_7+\sqrt{7}{C_{u}^{p}}^2 UU^{\top} u_7^{2} - 7 C_{u}^{p} R U u_7^{3}\right)}{T\left(C_{u}^{v} U + 7 R u_7^{2}\right)\left(\sqrt{7} C_{u}^{v}U - 7 C_{u}^{p} U u_7 + 7 \sqrt{7} R u_7^{2}\right)}\\~\\
	
	\frac{7 u_5\left(-C_{u}^{v} C_{u}^{p} UU^{\top} u_7+\sqrt{7}{C_{u}^{p}}^2 UU^{\top} u_7^{2} - 7 C_{u}^{p} R U u_7^{3}\right)}{T\left(C_{u}^{v} U + 7 R u_7^{2}\right)\left(\sqrt{7} C_{u}^{v}U - 7 C_{u}^{p} U u_7 + 7 \sqrt{7} R u_7^{2}\right)}\\~\\
	
	\frac{7 u_6\left(-C_{u}^{v} C_{u}^{p} UU^{\top} u_7+\sqrt{7}{C_{u}^{p}}^2 UU^{\top} u_7^{2} - 7 C_{u}^{p} R U u_7^{3}\right)}{T\left(C_{u}^{v} U + 7 R u_7^{2}\right)\left(\sqrt{7} C_{u}^{v}U - 7 C_{u}^{p} U u_7 + 7 \sqrt{7} R u_7^{2}\right)}\\~\\
	
	\frac{\sqrt{7} C_{u}^{p}  U u_7^{2}}{T\left( C_{u}^{v} U+7 R u_7^{2}\right)} \\~\\
	
	1
	\end{array}\right)$$
We can simply determine all of the system's eigenvalues and use the eigenvalues to get the system's left eigenvectors. We then gather it into the ODE system specified in Section 6, i.e.

\[
\begin{aligned}
	\frac{\mathrm{d} x}{\mathrm{~d} t} &=\frac{1}{\sigma_{k}(x, t, u)} \\
	\boldsymbol{l}_{k}(x, t, u)^{\top} A(x, t, u) \frac{\mathrm{d} u}{\mathrm{~d} t} &=\boldsymbol{l}_{k}(x, t, u)^{\top} g(x, t, u) \\
	k=1,2, \ldots, m .
\end{aligned}
\label{6} \tag{6}\]

\justify
Note that we have a point of singularity due to the fact that $\sigma_7 = 0$. We discuss how we deal with this in Section 11 and how the solution changes as we change the value of $\epsilon$. The above-reformulated system is easier to deal with numerically just for the fact that both $A$ and $B$ are almost in diagonal forms and that the eigenvalues and left-generalized eigenvectors are. \\

\justify
We could also proceed with the original system ie $M f_{t}+N f_{z}=g(t, z, f)$. First, we try to find the generalized eigenvalues and left eigenvectors of this system $\boldsymbol{l}_{k}^{\top}\left(M-\sigma_{k} N\right)=0,\; \forall k=1,2, \ldots, m$. In particular, the analytical solution to $\operatorname{det}(M-\sigma N)=0$ is done as follows

$$\begin{gathered}
	\operatorname{det}\left(\left(\begin{array}{cc}
		I-\frac{1}{\|F\|} F e^{\top} & -\frac{1}{T} F \\
		0^{\top} & \frac{c_{v}}{R T\|F\|} F
	\end{array}\right) - \sigma\left(\begin{array}{cc}
		\frac{R T}{A p} I & 0 \\
		\frac{R T}{A p} e^{\top} & \frac{c_{p}}{A p} F
	\end{array}\right)\right) = 0
\end{gathered}$$

$$\begin{gathered}
	= \operatorname{det}\left(\left(\begin{array}{cc}
		I-\frac{1}{\|F\|} F e^{\top} -\sigma \frac{R T}{A p} I & -\frac{1}{T} F \\
		-\sigma \frac{R T}{A p} e^{\top} & \frac{c_{v}}{R T\|F\|} F -\sigma \frac{c_{p}}{A p} F
	\end{array}\right)\right) = 0
\end{gathered}$$

\justify
Now using the properties of the block matrix we get that 

$$\begin{gathered}
 =  \operatorname{det}\left(\left(\begin{array}{cc}
	\left(I-\frac{1}{\|F\|} F e^{\top} -\sigma \frac{R T}{A p} I\right)\left(\frac{c_{v}}{R T\|F\|} F -\sigma \frac{c_{p}}{A p} F\right) - \left(-\frac{1}{T} F\right)\left(-\sigma \frac{R T}{A p} e^{\top}\right)
\end{array}\right)\right) = 0
\end{gathered}$$

\justify
From this easily we can figure out all the eigenvalues of the system and use the eigenvalues to figure out the left eigenvectors for this system. We then collect it into a system of ODEs as mentioned in Section 6 ie 

\[
\begin{aligned}
	\frac{\mathrm{d} z}{\mathrm{~d} t} &=\frac{1}{\sigma_{k}(z, t, f)} \\
	\boldsymbol{l}_{k}(z, t, f)^{\top} M(z, t, f) \frac{\mathrm{d} f}{\mathrm{~d} t} &=\boldsymbol{l}_{k}(z, t, f)^{\top} g(z, t, f) \\
	k=1,2, \ldots, m .
\end{aligned}
\label{7} \tag{7} \]
Both systems can, therefore then easily be solved by any numerical method that exists out there for ODEs or even PDEs. An attempt is also made to try to figure out a closed-form expression for a much simpler subset of the ODE system.

\subsection{Ordinary Differential Equation Solvers}

In the following section, we discuss some of the common techniques to solve the above system of ODE. We will proceed with trying to solve \ref{6} numerically. 

\subsection{Euler Methods}

We start off with one of the most common approaches to solving ODEs. The Euler method is a first-order numerical algorithm for solving ordinary differential equations (ODEs) using a known starting point. As described in \cite{11}, suppose we have the initial value problem to be solved on an interval $[t_0,t_f]$
$$
y^{\prime}(t)=f(t, y(t)), \quad y\left(t_{0}\right)=y_{0}
$$
Now divide this interval by the mesh points such that the value for $h = \frac{(t_f - t_0)}{N}$ is the size of every step and set $t_{n}=t_{0}+n h, \quad \forall n \in \{0,1,2,\ldots N\}$. Now, one step of the Euler method from $t_{n}$ to $t_{n+1}=t_{n}+h$ is
$$
y_{n+1}=y_{n}+h f\left(t_{n}, y_{n}\right), \quad \forall n \in \{0,\dots N-1\}
$$
The value of $y_{n}$ is an approximation of the solution to the ODE at time $t_{n}: y_{n} \approx y\left(t_{n}\right)$. The Euler method is explicit, i.e. the solution $y_{n+1}$ is an explicit function of $y_{i}$ for $i \leq n$ and this is known as the Euler's method. We may now get an approximation of a nearby point on a curve by going a short distance down a line tangent to the curve from any point on the curve. We begin by replacing the derivative in the differential equation above $y^{\prime}$ by the finite difference approximation and then integrating the differential equation between two consecutive mesh points $t_n$ and $t_{n+1}$
$$
y\left(t_{n+1}\right)=y\left(t_{n}\right)+\int_{t_{n}}^{t_{n+1}} f(t, y(t)) \mathrm{d} t, \quad n=0, \ldots, N-1,
$$
and then applying the numerical integration rule
$$
\int_{t_{n}}^{t_{n+1}} g(t) \mathrm{d} t \approx h g\left(t_{n}\right),
$$
with $g(t)=f(t, y(t))$, to get
$$
y\left(t_{n+1}\right) \approx y\left(t_{n}\right)+h f\left(t_{n}, y\left(t_{n}\right)\right), \quad n=0, \ldots N-1, \quad y\left(t_{0}\right)=y_{0} .
$$
To solve the mathematical model presented in the previous section ie section 6, a discretization of the space in the axial direction has been performed, as represented in the following equations \cite{6}:
$$
\begin{aligned}
	\frac{\partial u}{\partial x} &=\frac{U(x)-U(x-1)}{\Delta x} \\
	\frac{\partial T}{\partial x} &=\frac{T(x)-T(x-1)}{\Delta x}
\end{aligned}
$$
where $u, T$ are the molar concentration and the temperature, respectively. Specifically, the discretization process results in a 320 order system (seven chemical components and the temperature per 40 differential volumes) which has been implemented Mathematica $^{(\circledR)}$ and MATLAB $^{(\circledR)}$ and $\operatorname{SIMULINK}^{(B)}$ and integrated with the ode15s (stiff/NDF) solver and NDSolve feature. An example on how to use Mathematica's NDSolve method or even an implementation of the above method is as follows

\begin{lstlisting}
	(* steps = N *)
	(* initial = initialvalue *)
	myEuler[t0_, t1_, initial_, steps_] :=  
	With[{h = (t1 - t0)/steps},   FoldList[#1 + h f[#2, #1] &, initial, Range[t0 + h, t1, h]]]
	
	(* or we can use the NDSolve method*)
	NDSolve[\{y'[t] == y[t], x[0] == 1\}, x, \{t, 0, 10\}, StartingStepSize -> 1, Method -> ``ExplicitEuler'']
\end{lstlisting}

\subsection{Runge-Kutta}

We take a look at a generalization of the Euler technique. One of the most popular approaches for solving ODEs is the Runge-Kutta method. This method is a higher-order approximation of the halfway approach. Instead of shooting to the midpoint, computing the derivative, and then firing across the full interval, the Runge-Kutta technique takes four steps: shooting across one-quarter of the interval, estimating the derivative, shooting to the midway, and so on. To summarize, the goal of Runge–Kutta techniques is to approximate a Taylor series by executing successive (weighted) Euler steps. In this case, function evaluations (rather than derivatives) are used. This is simple to implement in Mathematica, as seen below.

\begin{lstlisting}
	NDSolve[{y'[t] == -y[t], y[0] == 1}, y[t], {t, 0, 1},
	Method -> "ExplicitMidpoint", "StartingStepSize" -> 1/10]
\end{lstlisting}
We can also extend this to incorporate a higher-order method by repeated function evaluation. 

\subsubsection{Remarks}

As a result, these approaches may be utilized to solve this system of $2m$ equations. The following equations are solved with an initialization of the state at the non-time coordinate's inlet boundary condition, and the results are gathered to generate ordered triplets $(x, t, u)$, which are vector $u$ values for each ordered pair $(x, t)$. The beginning values may be determined by selecting suitable values while keeping the molar and energy conservation equations in mind. In our case, we randomize the initial values to solve the above equation numerically. A more controlled approach would require the presence of a domain expert to let us know what the appropriate initial values for the system are.

\subsection{Partial Differential Equation Solvers}

In the following section, we discuss the most common technique to solve PDEs

\subsection{Finite Element Method}
Explicit closed-form solutions for partial differential equations (PDEs) are rarely available. The finite element method (FEM) is a technique to solve partial differential equations numerically. 

\justify
It is important for at least two reasons. First, the FEM is able to solve PDEs on almost any arbitrarily shaped region. Second, the method is well suited for use on a large class of PDEs. While it is almost always possible to conceive better methods for a specific PDE on a specific region, the finite element method performs quite well for a large class of PDEs.

\justify
It should be noted that the FEM can only provide an approximate answer. As a result, it is not the preferred method of resolving a physical condition. A closed-form analytical solution is the best technique to address a physical issue controlled by a differential equation. Unfortunately, in many real circumstances, the analytical answer is impossible to find or does not exist. This method can also be easily implemented in Mathematica as follows

\begin{lstlisting}
	Needs["NDSolve`FEM`"]
	(*Load the package FEM and try to solve the laplacian with boundary conditions*)
	pde = {Laplacian[u[x, y, z], {x, y, z}] == 1, 
	DirichletCondition[u[x, y, z] == 0, z == -1], 
	DirichletCondition[u[x, y, z] == 1, z == 1]};

	(*Solve the PDE using FEM*)
	AbsoluteTiming[
	MaxMemoryUsed[
	ILUOSolution = 
	NDSolveValue[pde, u, Element[{x, y, z}, mesh], 
	Method -> {"PDEDiscretization" -> {"FiniteElement", 
			"PDESolveOptions" -> {"LinearSolver" -> {Automatic, 
					Method -> {"Krylov", Method -> "BiCGSTAB", 
						Tolerance -> 10^-3, "Preconditioner" -> "ILU0"}}}}}]]/1024.^2]
\end{lstlisting}
\subsection{Tools}

\subsubsection{Mathematica}
The primary programming used for all of the analysis here is Mathematica. Mathematica is a very powerful language that can solve even handle symbolic calculations. More details of the various methods can be found in \cite{15}

\subsubsection{Matlab}

Matlab has a variety of methods that can solve not only systems of ODE but also systems of PDE. The primary two are the built-in functions ode23 and ode45, which implement Runge–Kutta 2nd/3rd-order and Runge–Kutta 4th/5th-order variants, respectively. The methods described are found in Matlab's documentation. The Matlab results were consistent with the results using Mathematica.

\section{Discussion}

We now proceed to numerically solve our reformulated system. We first consider the simplified version of the system  thereby considering only 2 molar equations ie Molar Flow Rate of $\mathrm{C_2H_5OH}$ and $\mathrm{CH_4}$ and the same temperature equation. Therefore we have our simplified system as $U = (u_1,u_2)$ and $u = (u_1,u_2,T)$
$$
\begin{gathered}
	A=\left(\begin{array}{ccc}
		1 & -\frac{u_{1}}{u_{2}} & -\frac{u_{1}}{T} \\
		0 & 0 & -\frac{u_{2}}{T} \\
		0 & 0 & \frac{C_{v}^{u} U}{\sqrt{2} R T u_{2}}
	\end{array}\right)
	B=\left(\begin{array}{ccc}
		\frac{RT}{Ap} & 0 & 0  \\
		0 & \frac{RT}{Ap} & 0 \\
		0 & \frac{\sqrt 2RT}{Ap} & \frac{C_{v}^{u} U}{Ap}
	\end{array}\right) \\
	g =\left(\begin{array}{cc}
		R T u_1 U V_1R' \\
		R T u_2 U V_2R' \\
		U \beta(T_f - T) - \sum_i u_i [R'H' - V_iR']^{T}
	\end{array}\right)
\end{gathered}
$$	
Again before attempting to solve the system we must first try to find the generalized eigenvalues and left eigenvectors so that we could transform our system of PDE into a system of ODEs. The eigenvalues of the above system are 

$$
\sigma_{1}=0, \quad \sigma{2} = \frac{Ap}{RT}, \quad  \sigma_{3}=\frac{A p}{R T C_{p}^{u} U}\left(\frac{C_{v}^{u} U}{\sqrt{2} u_{2}}-\sqrt{2} R u_{2}\right)
$$
\justify
and 
$$
\begin{gathered}
	l_{1}=\left(\begin{array}{c}
		\frac{u_1}{u_2} \\
		1 \\
		0
	\end{array}\right), l_{2}=\left(\begin{array}{c}
		1 \\
		0 \\
		0
	\end{array}\right), l_{1}=\left(\begin{array}{c}
	\frac{2 u_2\left(-C_{u}^{v} C_{u}^{p} UU^{\top} u_2+\sqrt{2}{C_{u}^{p}}^2 UU^{\top} u_2^{2} - 2 C_{u}
		^{p} R U u_2^{3}\right)}{T\left(C_{u}^{v} U + 2 R u_2^{2}\right)\left(\sqrt{2} C_{u}^{v}U - 2 C_{u}^{p} U u_2 + 2 \sqrt{2} R u_2^{2}\right)} \\
	 \frac{\sqrt{2} C_{u}^{p}  U u_2^{2}}{T\left( C_{u}^{v} U+2 R u_2^{2}\right)} \\
	1
\end{array}\right)
\end{gathered}
$$

\begin{figure}
	\centering
	\includegraphics[width=1\linewidth]{"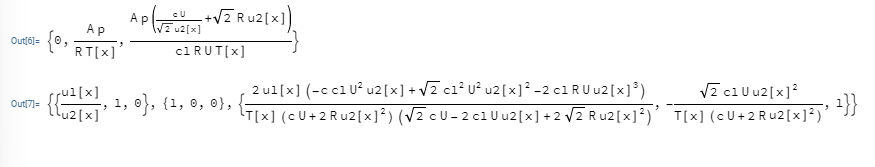"}
	\caption{Eigenvalues and Eigenvectors of our system}
	\label{fig:123}
\end{figure}

\justify
This is also confirmed in \ref{fig:123} where the constants $c$ and $c1$ are $C_{u}^{v} = c$ and $C_{u}^{p} = c1$. Now that we have our eigenvalues and eigenvectors we have the following ODE system
\[
\begin{aligned}
	\frac{\mathrm{d} x}{\mathrm{~d} t} &=\frac{1}{\sigma_{k}(x, t, u)} \\
	\boldsymbol{l}_{k}(x, t, u)^{\top} A(x, t, u) \frac{\mathrm{d} u}{\mathrm{~d} t} &=\boldsymbol{l}_{k}(x, t, u)^{\top} g(x, t, u) \\
	k=1,2,3.
\end{aligned}
\label{6} \tag{6}\]
\justify
which when explicitly written out can be written as

\[
\begin{aligned}
	\frac{\mathrm{d} x_3}{\mathrm{~d} t} &=\frac{1}{\sigma_{1}(x, t, u)} \\
	\frac{\mathrm{d} x_2}{\mathrm{~d} t} &=\frac{1}{\sigma_{3}(x, t, u)} \\
	\frac{\mathrm{d} x_1}{\mathrm{~d} t} &=\frac{1}{\sigma_{2}(x, t, u)} \\
	\boldsymbol{l}_{1}(x, t, u)^{\top} A(x, t, u) \frac{\mathrm{d} u_1}{\mathrm{~d} t} &=\boldsymbol{l}_{1}(x, t, u)^{\top} g(x, t, u) \\
	\boldsymbol{l}_{2}(x, t, u)^{\top} A(x, t, u) \frac{\mathrm{d} u_2}{\mathrm{~d} t} &=\boldsymbol{l}_{2}(x, t, u)^{\top} g(x, t, u) \\
	\boldsymbol{l}_{3}(x, t, u)^{\top} A(x, t, u) \frac{\mathrm{d} T}{\mathrm{~d} t} &=\boldsymbol{l}_{3}(x, t, u)^{\top} g(x, t, u) 
\end{aligned}
\label{6} \tag{6}\]

\justify
Notice that we have $\sigma_{1} = 0$ which implies that $\frac{\mathrm{d} x_1}{\mathrm{~d} t} =\frac{1}{\sigma_{1}(x, t, u)}$ would blow up and is a cause of singularity in our system. We replace the $0$ by $\epsilon > 0$ just to notice how the solution changes as $\epsilon$ increases. Finally, we make the following set of assumptions

\begin{enumerate}
	\item We now proceed to put the values for the constants ($R$, $T_f$, $\beta$, stoichiometric coefficients and etc) where some of the constants have the true values and others have been randomized. 
	\item We should also note that one of the eigenvalues in our system is $0$ which is a point of singularity. This was replaced by $\epsilon = 0.00001$ just to analyze the system better.
	\item Finally we set all the boundary conditions and initial conditions to small non-zero values but randomized.
	\item We set the independent variable $x$ to be solved in the range $[\epsilon,10]$
\end{enumerate}
\justify
We now use Mathematica's NDSolve feature to numerically solve the system by using the following code

\begin{lstlisting}
	(*Equations are defined in the appendix*)
	(*Initial conditions are all set to epsilon or 1*)
	(*Explictly without mentioning the Runge - Kutta Method was used*)
	eqns = {eqns1,eqns2,eqns3,eqns4,eqns5,eqns6,y1[0] == eps,y2[0] == 1,y3[0] == eps,u1[0] == eps,u2[0] == eps,T[0] == eps};
	sol = NDSolve[eqns,{y1,y2,y3,u1,u2,T},{x,eps,10}];
\end{lstlisting}

\justify
The above numerical method took exactly 113 steps to evaluate and converge to a solution. This was calculated as follows

\begin{lstlisting}
	(*Counter to keep track of the number of steps to converge*)
	Module[{c = 0}, 
	NDSolve[qns,{y1,y2,y3,u1,u2,T},{x,eps,10}, 
	EvaluationMonitor :> c++]; c]
\end{lstlisting}

\begin{figure}
	\centering
	\includegraphics[width=0.7\linewidth]{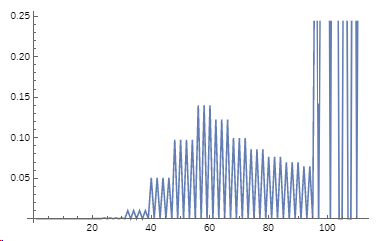}
	\caption{Successive differences vs steps}
	\label{fig:124}
\end{figure}

\justify
We also notice that during the evaluation of the numerical method, none of the evaluations failed at each step as noted in \ref{fig:124}.

\begin{figure}%
    \centering
    \subfloat[\centering $u_2$ solution for $x>0$]{{\includegraphics[width=8cm]{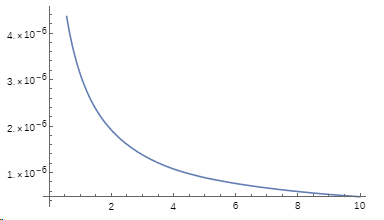} }}%
    \qquad
    \subfloat[\centering $u_2$ solution for $x<0$]{{\includegraphics[width=8cm]{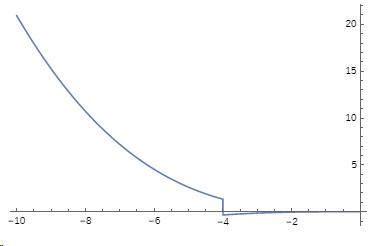} }}%
    \quad
    \subfloat[\centering $T$ solution for $x>0$]{{\includegraphics[width=8cm]{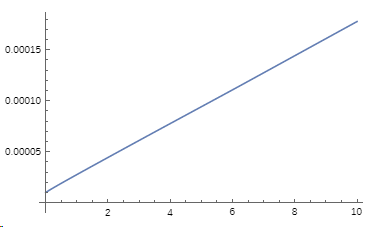} }}%
    \qquad
    \subfloat[\centering $T$ solution for $x<0$]{{\includegraphics[width=8cm]{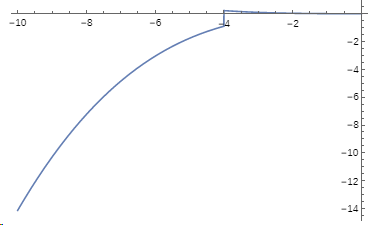} }}%
    \caption{Solution curves for $u_2$ and $T$.}%
    \label{fig:example2}%
\end{figure}

\begin{figure}%
    \centering
    \subfloat[\centering $x_1$ solution for $x>0$]{{\includegraphics[width=8cm]{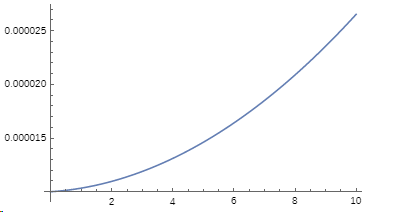} }}%
    \qquad
    \subfloat[\centering $x_1$ solution for $x<0$]{{\includegraphics[width=8cm]{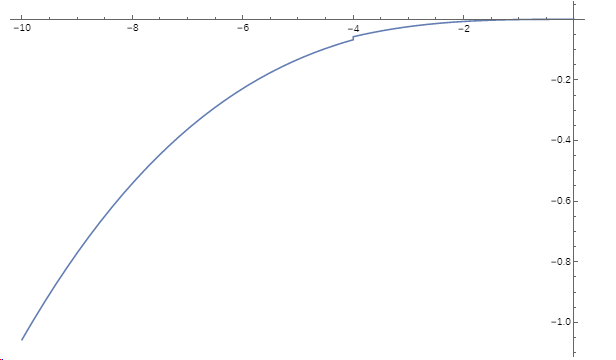} }}%
    \quad
    \subfloat[\centering $x_2$ solution for $x>0$]{{\includegraphics[width=8cm]{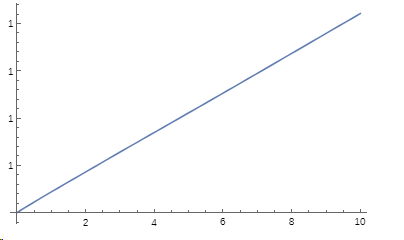} }}%
    \qquad
    \subfloat[\centering $x_2$ solution for $x<0$]{{\includegraphics[width=8cm]{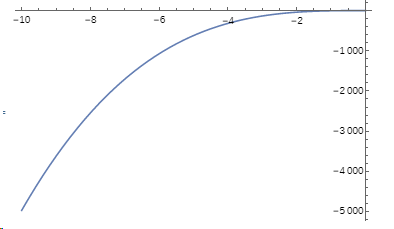} }}%
    \quad
    \subfloat[\centering $x_3$ solution for $x>0$]{{\includegraphics[width=8cm]{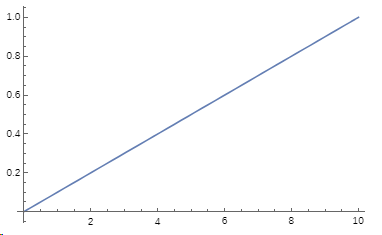} }}%
    \qquad
    \subfloat[\centering $x_3$ solution for $x<0$]{{\includegraphics[width=8cm]{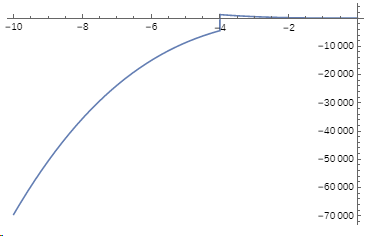} }}%
    \quad
    \subfloat[\centering $u_1$ solution for $x>0$]{{\includegraphics[width=8cm]{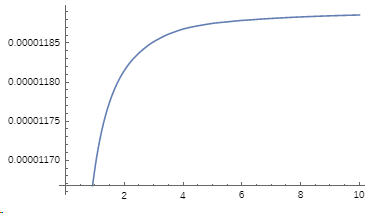} }}%
    \qquad
    \subfloat[\centering $u_1$ solution for $x<0$]{{\includegraphics[width=8cm]{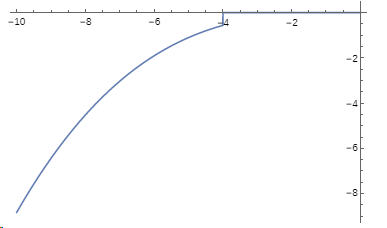} }}%
    \caption{Solution curves for $u_1$, $x_1$, $x_2$ and $x_3$.}%
    \label{fig:example}%
\end{figure}

\begin{figure}%
    \centering
    \subfloat[\centering All solutions plotted $x>0$ and $\epsilon = 0.00001$]{{\includegraphics[width=8cm]{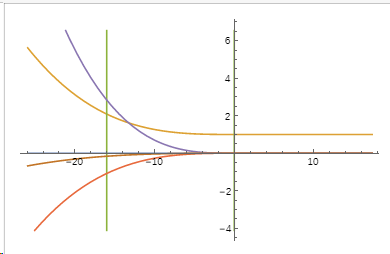} }}%
    \qquad
    \subfloat[\centering All solutions plotted $x<0$ and $\epsilon = 0.00001$]{{\includegraphics[width=8cm]{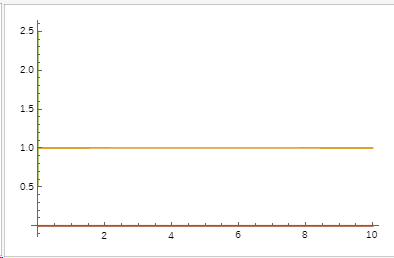} }}%
    \quad
    \subfloat[\centering All solutions plotted $x>0$ and $\epsilon = 1$]{{\includegraphics[width=8cm]{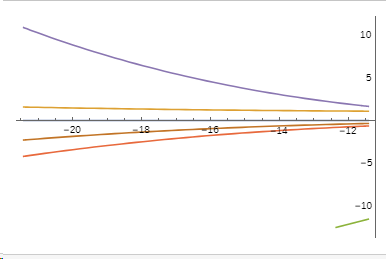} }}%
    \qquad
    \subfloat[\centering All solutions plotted $x<0$ and $\epsilon = 1$]{{\includegraphics[width=8cm]{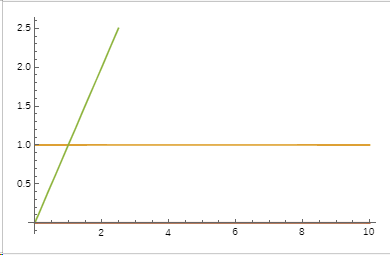} }}%
    \quad
    \subfloat[\centering All solutions plotted $x>0$ and $\epsilon = 100$]{{\includegraphics[width=8cm]{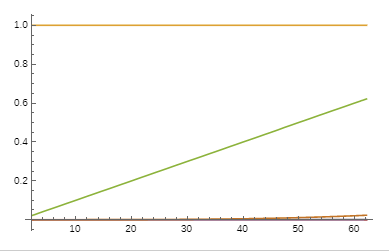} }}%
    \qquad
    \subfloat[\centering All solutions plotted $x<0$ and $\epsilon = 100$]{{\includegraphics[width=8cm]{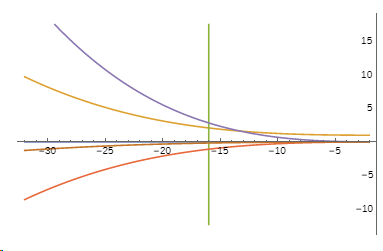} }}%
    \caption{All solution curves.}%
    \label{fig:example3}%
\end{figure}

\justify
Finally we plot all the solutions $x_1,x_2,x_3,u_1,u_2,T$ both for $x>0$ and $x<0$ as shown in figures \ref{fig:example2}, \ref{fig:example}, \ref{fig:example3}. Both the ODE Solvers Runge-Kutta and Euler methods result in the same numeric solutions and almost the same figures. We also consider changing the value of $\epsilon$ to see how the solution changes (all solutions plotted ie $x_1,x_2,x_3,u_1,u_2, T$ represented by the blue, yellow, green, red, purple, and orange lines respectively. We now note the following observations from the figures 

\begin{figure}
	\centering
	\includegraphics[width=1\linewidth]{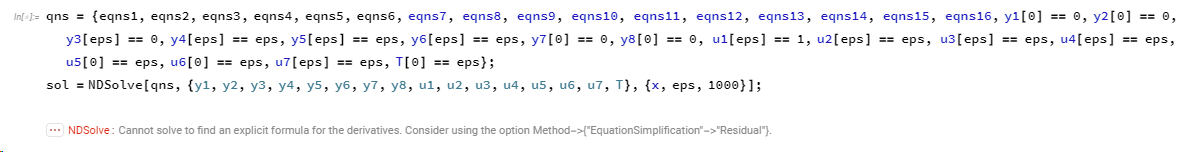}
	\caption{DSolve method to find the closed form solution}
	\label{fig:1j1}
\end{figure}

\begin{enumerate}
	\item As $\epsilon$ approaches $0$, the solution line $x_3$ approaches the $x$ axis.
	\item All the other solutions relatively remain the same and steady.
	\item When extrapolated into the negative $x$ axis we see that $x_3$ is highly random and discontinuous compared to the other solutions
	\item The solution has no closed-form solution. This is also confirmed in \ref{fig:1j1}.
	\item Increasing the number of evaluation steps and working precision gives a more accurate solution.
	\item We cannot gain a grasp of the numerical error of our system as there is no closed-form solution to compare it to, but we do know the truncation error (error made in a single evaluation) for Euler's Method is $O(h^2)$ and for Runge-Kutta Method is $O(h^p)$ where $p$ is the order of the method chosen.
	\item Overall the solution of the system is stable and this can be fed now into a Non-Linear Model Predictive control algorithm.
\end{enumerate}

\section{Limitations}
There were certain limitations and assumptions that were made due to lack of time; some of them were
\begin{itemize}
	\item The whole 16 equations were programmed, but we need a more detailed analysis and possibly a simplification of the equations for Mathematica even to find a numerical solution.
	\item Mathematica is similarly unable to handle an 8 by 8 system of PDEs. Therefore, we utilize the characteristics approach to convert it back to ODEs.
	\item Our reformulated system has a singularity point due to the fact that $\operatorname{det}(A) = 0$, and so we need to do more theoretical analysis on how it behaves.
	\item All of the constants were randomized; a more regulated setup is required to set out the constants as precisely as feasible in relation to the reactants.
\end{itemize}

\section{Future Work}

Although we did get some solid results in this project over the past semester, some of the future work for this research would include

\begin{itemize}
	\item An in-depth study of the system (more theoretical Analysis). This would involve proving the uniqueness theorem for the 7 by 7 system and how the $T$ equation behaves as $t \rightarrow \infty$. Hypothetical shock developments are also crucial to explore, as well as more attempts to convert the whole system into a conservation law.
	\item More research on the various numerical methods can help us to solve this system faster, as well as more numerical analysis.
	\item This technique may be generalized by combining these ODE solvers into Nonlinear Model Predictive Control algorithms that can be utilized for various Nonlinear systems.
	\item In this paper, we provided a static analysis of our model and, in the future, hope to provide a dynamic analysis of the model ie( when the time changes, how our solution evolves).
	\item Lastly, to completely numerically solve the 16 equations (Although I predict that it would just be the same as our simplified system with $u_7$ being our source of singularity and all the other solutions staying the same).
\end{itemize}

\section{Conclusion}

This paper provides a detailed analysis of the reduced-order model for an ethanol steam reformer that employs no approximations when converting a system of singular nonlinear distributed parameter systems to a system of nonlinear ordinary differential equations for computing the states required in nonlinear model predictive control formulations. We also go over the theoretical analysis of both the original and reformed systems, which is accomplished by establishing the uniqueness and existence theorem for quasi-linear and linear first-order systems. Despite the fact that we only have one point of singularity, we explain how to deal with it numerically. We also examined the potential of attempting to turn the system into a conservation law, which appeared to be impossible, but an attempt to convert the system into a conservation/convection law was successful. Finally, we offer some numerical analysis of the system and some techniques for solving it. Our results are consistent with those of prior works, allowing us to feed these ODE solvers into a Nonlinear Model Predictive Control method.

\justify
Because of the cheap online compute cost of the model equations (about $2m$ for ordinary differential equations), a mechanistic model may be employed in real-time control computations while explicitly accounting for input, state, and output limits. In conclusion, this research might be expanded to generalize the approach so that we could include additional reactor systems that are often used in chemical process control applications, for example. Finally, this would move us closer to our aim of producing hydrogen safely, which might then be utilized as green energy.

\section*{Acknowledgements}~\\
I thank Professor Xinwei Yu for supervising me this semester; without him, I wouldn't have been able to do this. I also want to thank Professor Paul Buckingham for all of his assistance and criticism during the semester, as well as for the fantastic course. I'd also like to thank all of the lecturers that came to class to deliver presentations, as well as my peers and family for supporting me throughout this semester.

\justify
Finally, I just wanted to say thank you to everyone who is a part of my life and who has stuck by my side and helped me when things became tough this semester. I thank God for keeping me alive and for providing me with the strength and wisdom I need to live my life in a fruitful way.

~\\

\appendix

\pagebreak
\section*{Appendices}
\addcontentsline{toc}{section}{Appendices}
\renewcommand{\thesubsection}{\Alph{subsection}}

\subsection{Eigenvalues and Eigenvectors}
\begin{lstlisting}
	(*Define the matrices A,B,g and the constants values have been randomized*)
	A = 100;
	V1 =  V2 =  V3 = v4 = V5 = V6 = V7  =1;
	c = c1 = 100;
	R' = R = 10;
	H = 10;
	A =\begin{center}
		\begin{flushleft}
			
		\end{flushleft}
	\end{center} {{1,0,0,0,0,0,-u1[x]/u7[x],-u1[x]/T[x]},{0,1,0,0,0,0,-u2[x]/u7[x],-u2[x]/T[x]},{0,0,1,0,
		0,0,-u3[x]/u7[x],-u3[x]/T[x]},{0,0,0,1,0,0,-u4[x]/u7[x],-u4[x]/T[x]},{0,0,0,0,1,0,
		-u5[x]/u7[x],-u5[x]/T[x]},{0,0,0, 0,0,1,-u6[x]/u7[x],-u6[x]/T[x]},{0,0,0,0,0
		,0,0,-u7[x]/T[x]},{0,0,0,0,0,0,0,(c U)/(Sqrt[7] R T[x] u7[x])}}
	
	B = {{(R T[x])/(A p), 0,0,0,0,0,0,0},{0,(R T[x])/(A p), 0,0,0,0,0,0}, {0,0, (R T[x])/(A p), 0,0,0,0,0}, {0,0,0,(R T[x])/(A p), 0,0,0,0}, {0,0,0,0,(R T[x])/(A p), 0,0,0}, {0,0,0,0,0,(R T[x])/(A p), 0,0}, {0,0,0,0,0,0,(R T[x])/(A p), 0},{0,0,0,0,0,0,(Sqrt[7] R T[x])/(A p),(c1 U)/(A p)}}
	
	
	g = {{R T[x] (u1[x] u1[x] + u1[x] u2[x] +u1[x] u3[x] + u1[x] u4[x] + u1[x] u5[x] + u1[x] u6[x] +u1[x] u7[x]) V1 R'},{R T[x](u2[x] u1[x] + u2[x] u2[x] +u2[x] u3[x] + u2[x] u4[x] + u2[x] u5[x] + u2[x] u6[x] +u2[x] u7[x]) V2 R'},{R T[x] (u3[x] u1[x] + u3[x] u2[x] +u3[x] u3[x] + u3[x] u4[x] + u3[x] u5[x] + u3[x] u6[x] +u3[x] u7[x]) V3 R'},{R T[x] (u4[x] u1[x] + u4[x] u2[x] +u4[x] u3[x] + u4[x] u4[x] + u4[x] u5[x] + u4[x] u6[x] +u4[x] u7[x]) V4 R'},{R T[x] (u5[x] u1[x] + u5[x] u2[x] +u5[x] u3[x] + u5[x] u4[x] + u5[x] u5[x] + u5[x] u6[x] +u5[x] u7[x]) V5 R'},{R T[x] (u6[x] u1[x] + u6[x] u2[x] +u6[x] u3[x] + u6[x] u4[x] + u6[x] u5[x] + u6[x] u6[x] +u6[x] u7[x]) V6 R'},{R T[x] (u7[x] u1[x] + u7[x] u2[x] +u7[x] u3[x] + u7[x] u4[x] + u7[x] u5[x] + u7[x] u6[x] +u7[x] u7[x]) V7 R'},{U1 B (Tf - T[x]) -u1[x] (R' H - V1 R T[x]) - u2[x](R' H - V2 R T[x]) - u3[x](R' H - V3 R T[x]) - u4[x](R' H - V4 R T[x]) - u5[x](R' H - V5 R T[x]) - u6[x](R' H - V6 R T[x]) - u7[x](R' H - V7 R T[x])}} 
	
	(*Find the generalized Eigenvalues and Eigenvectors of our system*)
	l = Eigenvectors[{a,b}]
	s = Eigenvalues[{a,b}]
\end{lstlisting}
\subsection{ODE Solver}
\begin{lstlisting}
	
	(*Define the 16 equations of our system *)
	eqns1 = y1'[x] == 1/(s[[2]]);
	eqns2 = y2'[x]  == 1/(s[[3]]);
	eqns3 = y3'[x]  == 1/(s[[4]]);
	eqns4 = y4'[x]  == 1/(s[[5]]);
	eqns5 = y5'[x]  == 1/(s[[6]]);
	eqns6 = y6'[x]  == 1/(s[[7]]);
	eqns7 = y7'[x]  == 1/(s[[1]] + 1);
	eqns8 = y8'[x] == 1/(s[[8]]);
	eqns9 = u1'[x] == (R T[x] (u1[x] u6[x] + u2[x] u6[x] +u3[x] u6[x] + u4[x] u6[x] + u5[x] u6[x] + u6[x] u6[x] +u7[x] u6[x]) V6 R');
	eqns10 = u2'[x] == (R T[x] (u1[x] u5[x] + u2[x] u5[x] +u5[x] u5[x] + u4[x] u5[x] + u5[x] u5[x] + u6[x] u5[x] +u7[x] u5[x]) V5 R');
	eqns11 = u3'[x] ==  (R T[x] (u1[x] u4[x] + u2[x] u4[x] +u5[x] u4[x] + u4[x] u4[x] + u5[x] u4[x] + u6[x] u4[x] +u7[x] u4[x]) V4 R');
	eqns12 = u4'[x] ==  (R T[x] (u1[x] u3[x] + u2[x] u3[x] +u5[x] u3[x] + u4[x] u3[x] + u5[x] u3[x] + u6[x] u3[x] +u7[x] u3[x]) V3 R');
	eqns13 = u5'[x] ==  (R T[x] (u1[x] u2[x] + u2[x] u2[x] +u5[x] u2[x] + u4[x] u2[x] + u5[x] u2[x] + u6[x] u2[x] +u7[x] u2[x]) V2 R');
	eqns14 = u6'[x] ==  (R T[x] (u1[x] u1[x] + u2[x] u1[x] +u5[x] u1[x] + u4[x] u1[x] + u5[x] u1[x] + u6[x] u1[x] +u7[x] u1[x]) V1 R');
	array1 = l[[8]]*g;
	ufinal = {u1'[x],u2'[x],u3'[x],u4'[x],u5'[x],u6'[x],u7'[x],T'[x]};
	new1 = ArrayReshape[array1, Dimensions[array1] ~DeleteCases~ 1];
	eqns15 = Total[l[[8]].a ufinal] == Total[new1];
	eqns16 = -((u1[x]u1[x]/(u7[x] u7[x])) +(u2[x]u2[x]/(u7[x] u7[x])) + (u3[x]u3[x]/(u7[x] u7[x])) + (u4[x]u4[x]/(u7[x] u7[x])) + (u5[x]u5[x]/(u7[x] u7[x])) + (u6[x]u6[x]/(u7[x] u7[x]))) u7'[x] == (((R T[x] u6[x] (u1[x] u6[x] + u2[x] u6[x] +u3[x] u6[x] + u4[x] u6[x] + u5[x] u6[x] + u6[x] u6[x] +u7[x] u6[x]) V6 R')/(u7[x])) + ((R T[x] u5[x] (u1[x] u5[x] + u2[x] u5[x] +u3[x] u5[x] + u5[x] u6[x] + u5[x] u5[x] + u6[x] u5[x] +u7[x] u5[x]) V5 R')/(u7[x])) + ((R T[x] u4[x] (u1[x] u4[x] + u2[x] u4[x] +u3[x] u4[x] + u4[x] u4[x] + u5[x] u4[x] + u6[x] u4[x] +u7[x] u4[x]) V4 R')/(u7[x])) + ((R T[x] u3[x] (u1[x] u3[x] + u2[x] u3[x] +u3[x] u3[x] + u4[x] u3[x] + u5[x] u3[x] + u6[x] u3[x] +u7[x] u3[x]) V3 R')/(u7[x])) + ((R T[x] u2[x] (u1[x] u2[x] + u2[x] u2[x] +u3[x] u2[x] + u4[x] u2[x] + u5[x] u2[x] + u6[x] u2[x] +u7[x] u2[x]) V2 R')/(u7[x])) + ((R T[x] u1[x] (u1[x] u1[x] + u2[x] u1[x] +u3[x] u1[x] + u1[x] u1[x] + u1[x] u6[x] + u1[x] u1[x] +u1[x] u6[x]) V1 R')/(u7[x])) + R T[x] (u1[x] u7[x]+u2[x] u7[x]+u3[x] u7[x]+u4[x] u7[x]+u5[x] u7[x]+u6[x] u7[x]+u7[x] u7[x]) V7 R');
	eps = 10^-5; 
	
	(*Solve the system numerically with initial conditions*)
	qns = {eqns1,eqns2,eqns3,eqns4,eqns5,eqns6,eqns7, eqns8, eqns9, eqns10, eqns11, eqns12, eqns13, eqns14, eqns15, eqns16,y1[0] == 0,y2[0] == 0,y3[eps] == 0,y4[eps] == eps,y5[eps] == eps,y6[eps] == eps,y7[0] == 0,y8[0] == 0, u1[eps] == 1,u2[eps] == eps, u3[eps] == eps,u4[eps] == eps, u5[0] == eps,u6[0] == eps,u7[eps] == eps,T[0] == eps};
	sol =NDSolve[qns,{y1,y2,y3,y4,y5,y6,y7,y8,u1,u2,u3,u4,u5,u6,u7,T},{x,eps,1000}];
	
\end{lstlisting}

\subsection{Simplified system}

\begin{lstlisting}
	(*Definining our system*)
	a1= {{1,-u1[x]/u2[x],-u1[x]/T[x]},{0,0 ,-u2[x]/T[x]},{0,0,(c U)/(Sqrt[2] R T[x] u2[x])}}
	
	b1= {{(R T[x])/(A p), 0,0},{0,(R T[x])/(A p), 0},{0,(Sqrt[2] R T[x])/(A p),(c1 U)/(A p)}}
	
	g1 =  {{R T[x] (u1[x] u1[x] + u1[x] u2[x] ) V1 R'},{R T[x] (u2[x] u1[x] + u2[x] u2[x] ) V2 R'},{U1 B (Tf - T[x]) -u1[x] (R' H - V1 R T[x]) - u2[x](R' H - V2 R T[x]) }}
	
	(*Eigenvalues and Eigenvectors of our system*)
	s1 = Eigenvalues[{a1,b1}]
	l1 = Eigenvectors[{a1,b1}]
	
	(*Define constants and our 6 equations*)
	array1 = l1[[3]]*g1;
	new1 = ArrayReshape[array1, Dimensions[array1] ~DeleteCases~ 1];
	new1;
	epsilon =1000000;
	R = 6.022*10;
	R' = 1;
	V1 = 2;
	eps = 10^-5;
	inf = 5;
	V2 = 3;
	c1 = 10;
	U =1 ;
	H = 1;
	Tf = 232;
	B = 10;
	U1 = 10;
	c = 10;
	A = 343;
	p = 10;
	u10 = {u1'[x],u2'[x],T'[x]};
	eqns1 = y1'[x] == 1/(s1[[2]]);
	eqns2 = y2'[x]  == 1/(s1[[3]]);
	eqns3 = y3'[x] == 1/(s1[[1]]+epsilon);
	eqns4 = u1'[x] + u1[x] (u2'[x]/u2[x]) + u1[x] (T'[x]/T[x]) == (R T[x] (u1[x]u1[x] + u1[x]u2[x]) V1 R');
	eqns5 =(-(u1[x]^2/(T[x] u2[x])) - u2[x]/T[x]) Derivative[1][T][x] + (u1[x] Derivative[1][u1][x])/u2[x] - (u1[x]^2 Derivative[1][u2][x])/u2[x]^2 == (R T[x] u1[x] (u1[x]^2 + u1[x] u2[x]) V1 Derivative[1][R])/u2[x] +  R T[x] (u1[x] u2[x] + u2[x]^2) V2 Derivative[1][R];
	eqns6 = Total[l1[[3]].a1 u10]  == Total[new1];
	
	(*Solution of our system*)
	qns = {eqns1,eqns2,eqns3,eqns4,eqns5,eqns6,y1[0] == eps,y2[0] == 1,y3[0] == eps,u1[0] == eps,u2[0] == eps,T[0] == eps};
	sol = NDSolve[qns,{y1,y2,y3,u1,u2,T},{x,eps,10},Method -> "ExplicitRungeKutta"];
	
	(*2D Plot of all our solutions*)
	DynamicModule[{xc,dx},Manipulate[Plot[Evaluate[{y1[x],y2[x],y3[x],u1[x],u2[x],T[x]}/. sol[[1]]],{x,xc-dx,xc+dx}],
	{{xc,1000001/200000,"center"},-(4999993/200000),1399999/40000},{{dx,5.,"zoom"},2999997/100000,0.3}]
	,DynamicModuleValues:>{}] 
\end{lstlisting}

\end{document}